\documentclass[10pt]{article}
\usepackage{amstext}
\usepackage{amsfonts}
\usepackage{amssymb}
\usepackage{amsbsy} 
\usepackage{latexsym}
\usepackage{xy}
\usepackage{hhline}
\xyoption{all}
\newcommand{\vtx}[1]{*+[o][F-]{\scriptscriptstyle #1}} 
\newcounter{num}[section] %

\newenvironment{theo}
{\refstepcounter{num}%
\bigskip\noindent{\bf Theorem~\arabic{section}.\arabic{num}. }\it}

\newenvironment{prop}
{\refstepcounter{num}%
\bigskip\noindent{\bf Proposition~\arabic{section}.\arabic{num}. }\it}

\newenvironment{cor}
{\refstepcounter{num}%
\bigskip\noindent{\bf Corollary~\arabic{section}.\arabic{num}. }\it}

\newenvironment{lemma}
{\refstepcounter{num}%
\bigskip\noindent{\bf Lemma~\arabic{section}.\arabic{num}. }\it}

\newcommand{\example}
{\refstepcounter{num}%
\bigskip\noindent{\bf Example~\arabic{section}.\arabic{num}.}}
\newcommand{\remark}
{\refstepcounter{num}%
\bigskip\noindent{\bf Remark~\arabic{section}.\arabic{num}.}}
\newcommand{\definition}[1]
{\refstepcounter{num}%
\bigskip\noindent{\bf Definition~\arabic{section}.\arabic{num}}~({\it #1}).}

\newenvironment{proof}{\noindent{\it Proof. }}
{$\Box$ \bigskip}

\newenvironment{eq}{\begin{equation}}{\end{equation}}

\newcommand{\Ref}[1]{(\ref{#1})}

\newcommand{\si}{\sigma}
\newcommand{\al}{\alpha}
\newcommand{\be}{\beta}
\newcommand{\ga}{\gamma}
\newcommand{\la}{\lambda}
\newcommand{\de}{\delta}

\newcommand{\LA}{\langle}
\newcommand{\RA}{\rangle}
\newcommand{\ov}[1]{\overline{#1}}
\newcommand{\un}[1]{{\underline{#1}} }

\newcommand{\tr}{\mathop{\rm tr}}

\newcommand{\mdeg}{\mathop{\rm mdeg}}
\newcommand{\diag}{\mathop{\rm diag}}

\newcommand{\sign}{\mathop{\rm{sgn }}}
\newcommand{\Hom}{{\mathop{\rm{Hom }}}}
\newcommand{\End}{{\mathop{\rm{End }}}}

\newcommand{\rep}[2]{{\mathop{\rm{rep}}_{#1}}(#2)}
\newcommand{\repO}[2]{{\mathop{\rm{rep}}_{#1}^{\rm O}}(#2)}
\newcommand{\repSp}[2]{{\mathop{\rm{rep}}_{#1}^{\rm Sp}}(#2)}

\newcommand{\algA}{\mathcal{A}}    
\newcommand{\algB}{\mathcal{B}}    

\newcommand{\M}{\mathcal{M}}  
\newcommand{\T}{{\mathcal T} } 

\newcommand{\FF}{{\mathbb{F}}}   
\newcommand{\NN}{{\mathbb{N}} }
\newcommand{\ZZ}{{\mathbb{Z}} }  
\newcommand{\QQ}{{\mathbb{Q}} }

\renewcommand{\P}{{\mathop{\ov{\rm{pf}}}}}
\newcommand{\Pf}{{\mathop{\rm{pf }}}}
\newcommand{\DP}{{\rm DP} }
\newcommand{\F}{{\mathop{\rm{bpf }}}}

\newcommand{\ovphi}[1]{\varphi(#1)}

\newcommand{\QS}{\mathfrak{Q}} 

\newcommand{\Q}{\mathcal{Q}}    
\newcommand{\n}{\boldsymbol{n}} 
\newcommand{\g}{\boldsymbol{g}} 
\newcommand{\h}{\boldsymbol{h}} 
\renewcommand{\i}{\boldsymbol{i}} 

\newcommand{\D}{{D}}

\newcommand{\Sp}{S\!p}

\newcommand{\Symmgr}{\mathcal{S}}                 

\newcommand{\loopR}[3]{%
\begin{picture}(20,0)(#1,#2)
\put(-2,1){\llap{$\scriptstyle #3$}} \put(11,3){\circle{20}} \put(20,6){\vector(1,-4){1}}
\end{picture}}
\newcommand{\loopL}[3]{%
\begin{picture}(20,0)(#1,#2)
\put(22,1){$\scriptstyle #3$} \put(9,3){\circle{20}} \put(0,6){\vector(-1,-4){1}}
\end{picture}}

\newcommand{\rectangle}[2]{
\begin{picture}(0,0)
\put(-#1,-#2){\line(1,0){#1}}\put(0,-#2){\line(1,0){#1}}
\put(-#1,#2){\line(1,0){#1}}\put(0,#2){\line(1,0){#1}}
\put(-#1,-#2){\line(0,1){#2}}\put(-#1,0){\line(0,1){#2}}
\put(#1,-#2){\line(0,1){#2}}\put(#1,0){\line(0,1){#2}}
\end{picture}}

\begin{document}
\title{Representations of quivers, their generalizations and invariants}
 \author{}
\date{} 
\maketitle

\vspace{-1.5cm} 
$$
\begin{array}{cc}
\mbox{ \large A. A. Lopatin}&\mbox{\large A. N. Zubkov} \\
\mbox{ Institute of Mathematics ,}&\mbox{ Omsk State Pedagogical University,}\\
\mbox{ Siberian Branch of }&\mbox{Chair of Geometry,}\\
\mbox{ the Russian Academy of Sciences, }&\mbox{}\\
\mbox{ Pevtsova street, 13,}&\mbox{Tuhachevskogo embarkment, 14,}\\
\mbox{ Omsk 644099 Russia}&\mbox{Omsk 644099, Russia} \\
\mbox{ artem\underline{ }lopatin@yahoo.com}&\mbox{a.zubkov@yahoo.com} \\
\mbox{http://www.iitam.omsk.net.ru/\~{}lopatin/}&\\
\end{array}
$$

\begin{abstract} This paper is a survey on invariants of representations of quivers and their
generalizations such as (super)mixed representations by Zubkov, orthogonal and symplectic representations by Derksen and Weyman, and representations of signed quivers by Shmelkin. Working over a field of arbitrary characteristic we present the description of generating systems for invariants and
relations between generators.
\end{abstract}



\tableofcontents

\section{Introduction}\label{section_intro}

We work over an infinite field $\FF$ of arbitrary characteristic. All vector spaces are
supposed to be finite dimensional and all algebras are associative with unity.

The paper can be split into three parts. The first part is Section~\ref{section_matrix}
and it is dedicated to matrix invariants of classical linear groups. In
Section~\ref{subsection_matrix} systems of generators for algebras of matrix invariants
are given.  Some connections between matrix invariants and representations of free
algebras are considered in Sections~\ref{subsection_alg} and~\ref{subsection_inv}.
Relations between generators for the cases of $GL(n)$ and $O(n)$ are described in
Section~\ref{subsection_relations}. The description in the case of $O(n)$ is given in
terms of the polynomial $\si_{t,r}$, which is defined in Section~\ref{subsection_DP}.

The second part is Section~\ref{section_quivers} and it is dedicated to representations
of quivers. A {\it quiver} is a finite oriented graph. A representation of dimension
$(\n_1,\ldots,\n_l)$ of a quiver with $l$ vertices consists of a collection of column
vector spaces $\FF^{\n_1},\ldots,\FF^{\n_l}$, assigned to the vertices, and linear
mappings between the vector spaces ``along"{} the arrows (see
Section~\ref{subsection_quivers_definitions} for details). This notion was introduced
by Gabriel in~\cite{Gabriel72} as an effective mean for description of some problems of
the linear algebra. The importance of this notion from the representation theory point of
view is due to the following fact. Let $\algA$ be a finite dimensional basic algebra over
algebraically closed field. Then the category of finite dimensional modules over $\algA$
is a full subcategory of the category of representations of some
quiver (see Chapter~3 from~\cite{Kirichenko}). 
Invariants of quivers are important not only in the invariant theory but also in the
representational theory of quivers because these invariants distinguish semi-simple
representations of a quiver. A system of generators for invariants of a quiver is given
in Section~\ref{subsection_quivers_generators} and relations between generators are
described in Section~\ref{subsection_quivers_relations}.

The third part consists of Sections~\ref{section_mixed_quivers},~\ref{section_bpf},
and~\ref{section_generators} and it is dedicated to generalizations of representations of
quivers. Given a representation of dimension $(\n_1,\ldots,\n_l)$ of a quiver with $l$
vertices, we generalize this notion as follows. Let $v$ be a vertex of the quiver, where
$1\leq v\leq l$. In the classical case $GL(\n_v)$ acts on $\FF^{\n_v}$ but in our case an
arbitrary classical linear group from the list $GL(\n_v)$, $O(\n_v)$, $Sp(\n_v)$,
$SL(\n_v)$, $SO(\n_v)$ can act on $\FF^{\n_v}$. Moreover, we can consider the dual space
$(\FF^{\n_v})^{\ast}$ together with $\FF^{\n_v}$ in order to deal with bilinear forms
together with linear mappings. Finally, instead of arbitrary linear mappings ``along"{}
arrows we can consider only those that, for example, preserve some bilinear symmetric
form on ``vertex"{} spaces, etc. Representations of quivers obtained in this way are
called {\it $\QS$-mixed representations}, where a {\it mixed quiver setting} $\QS$ is the
data that determines what linear mappings we may take to form a mixed representation,
i.e., $\QS$ is a quiver together with a dimension vector, a product $G$ of classical
linear groups, kinds of linear mappings along arrows and an involution that shows which
vertices are dual. The exact definitions together with examples are given in
Section~\ref{subsection_mixed_quivers_definitions}.

Particular cases of the considered construction are
\begin{enumerate}
\item[$\bullet$] {\it mixed} and {\it supermixed} representations of quivers, introduced
by Zubkov in~\cite{Zubkov_preprint} (see also~\cite{ZubkovI}); to obtain mixed
(supermixed, respectively) representations we should take $\QS$-mixed representations
such that $G$ is equal to a product of the general linear groups (the general linear
groups, orthogonal and symplectic groups, respectively);

\item[$\bullet$] {\it orthogonal} and {\it symplectic} representations of {\it symmetric}
quivers, introduced by Derksen and Weyman in~\cite{DW02};

\item[$\bullet$] representations of {\it signed} quivers, introduced by Shmelkin
in~\cite{Shmelkin};
\end{enumerate}
(see Example~\ref{ex3part4} for details). The motivation for these generalizations of
quivers from the point of view of the representational theory of algebraic groups was
given in~\cite{DW02},~\cite{Shmelkin}, where symmetric and signed quivers, respectively,
of tame and finite type were classified.

Let us depict different kinds of algebras of invariants and semi-invariants.  Here a vertical arrow means that the upper class of algebras contains the lower class of algebras. All necessary definitions can be found in Sections~\ref{subsection_matrix}, 
\ref{subsection_quivers_definitions}, and~\ref{subsection_mixed_quivers_definitions}
(see also Examples~\ref{ex3part4} and~\ref{ex_semi}). %
$$
\begin{picture}(0,250)
\put(0,230){%
\put(-95,0){
\put(0,0){\small Invariants of $\QS$-mixed representations of quiver}
\put(95,3){\rectangle{200}{10}}
} %
\put(-160,-35){
\put(0,0){\small Invariants of supermixed}
\put(-1,-10){\small representations of quivers}
\put(50,-2){\rectangle{90}{15}}
} %
\put(50,-35){
\put(0,0){\small Semi-invariants of supermixed}
\put(10,-10){\small representations of quivers}
\put(60,-2){\rectangle{90}{15}}
} %
\put(-80,-80){
\put(0,0){\small Invariants}
\put(3,-10){\small of mixed}
\put(-9,-20){\small representations}
\put(1,-30){\small of quivers}
\put(20,-12){\rectangle{40}{25}}
} %
\put(-80,-145){
\put(-4,0){\small Invariants of}
\put(-9,-10){\small representations}
\put(2,-20){\small of quivers}
\put(20,-7){\rectangle{40}{20}}
} %
\put(-80,-200){
\put(5,0){\small Matrix}
\put(-15,-10){\small $GL(n)$-invariants}
\put(20,-2){\rectangle{40}{15}}
} %
\put(-195,-200){
\put(0,0){\small Matrix $O(n)$- and}
\put(4,-10){\small $\Sp(n)$-invariants}
\put(35,-2){\rectangle{40}{15}}
} %
\put(140,-80){
\put(-10,0){\small Semi-invariants}
\put(3,-10){\small of mixed}
\put(-9,-20){\small representations}
\put(1,-30){\small of quivers}
\put(20,-12){\rectangle{40}{25}}
} %
\put(140,-145){
\put(-15,0){\small Semi-invariants of}
\put(-9,-10){\small representations}
\put(2,-20){\small of quivers}
\put(20,-7){\rectangle{40}{20}}
} %
\put(40,-200){
\put(5,0){\small Matrix}
\put(-15,-10){\small $SO(n)$-invariants}
\put(20,-2){\rectangle{40}{15}}
} %
\put(-110,-7){\vector(0,-1){15}}
\put(110,-7){\vector(0,-1){15}}
\put(-160,-52){\vector(0,-1){135}}
\put(-60,-52){\vector(0,-1){15}}
\put(60,-52){\vector(0,-1){135}}
\put(160,-52){\vector(0,-1){15}}
\put(-60,-117){\vector(0,-1){15}}
\put(160,-117){\vector(0,-1){15}}
\put(-60,-172){\vector(0,-1){15}}
}%
\end{picture}
$$

In Section~\ref{subsection_mixed_quivers_generators} generators for invariants of
supermixed representations of quivers are given. Relations between generators for mixed
representations of quivers are described in
Section~\ref{section_mixed_quivers_relations}. Earlier introduced polynomial $\si_{t,r}$ plays a key role here.

To describe generators for invariants of $\QS$-mixed representations, we introduce a
block partial linearization of the pfaffian (b.p.l.p.) and a tableau with substitution
(see Section~\ref{section_bpf}). Note that the polynomial $\si_{t,r}$ can be defined in
terms of a b.p.l.p.~(see part~4 of Example~\ref{ex_pf_1} and Remark~\ref{rem_history001}).

In Section~\ref{section_generators} we present generators for invariants of $\QS$-mixed
representations. Before formulating the main result in
Section~\ref{subsection_general_case}, we consider a partial case of semi-invariants for
the usual representations of bipartite quivers. Note that Theorem~\ref{theo_main} implies
the description of semi-invariants of supermixed representations of quivers.

So the main results formulated in this paper are the descriptions of
\begin{enumerate}
\item[$\bullet$] generators for invariants of $\QS$-mixed representations from Section~\ref{subsection_general_case};

\item[$\bullet$] relations for invariants of mixed representations of quivers from 
Section~\ref{section_mixed_quivers_relations};

\item[$\bullet$] relations for matrix $O(n)$-invariants from Section~\ref{subsection_relations}.
\end{enumerate} %
The rest of statements from this paper that concern generators or relations are partial cases of the mentioned ones.

\section{Notations}\label{section_notations}

\subsection{Matrices}\label{subsection_matrices}
Given a positive integer $n$, let us fix the following notations for the
classical linear groups:
\begin{enumerate}
\item[$\bullet$] $O(n)=\{A\in \FF^{n\times n}\,|\,A A^T=A^T A=E\}$ and $SO(n)=\{A\in
O(n)\,|\,\det(A)=1\}$, where we assume that the characteristic of $\FF$ is not two;

\item[$\bullet$] $\Sp(n)=\{A\in \FF^{n\times n}\,|\,A^T J A=J\}$, where we assume that
$n$ is even.
\end{enumerate}
Here $E=E(n)$ stands for the identity matrix, and $J=J(n)=\left(
\begin{array}{cc}
0& E(n/2) \\
-E(n/2)& 0\\
\end{array}
\right)$ stands for the the matrix of the skew-symmetric bilinear form on $\FF^{2n}$. We
also fix the following notations for certain subspaces of $\FF^{n\times n}$:
\begin{enumerate}
\item[$\bullet$] $S^{+}(n)=\{A\in \FF^{n\times n} \,|\, A^T=A\}$ is the space of
symmetric matrices;

\item[$\bullet$] $S^{-}(n)=\{A\in \FF^{n\times n} \,|\, A^T=-A\}$ is the space of
skew-symmetric matrices;

\item[$\bullet$] $L^{+}(n)= %
\{A\in \FF^{n\times n} \,|\, AJ\text{ is a symmetric matrix}\}$;

\item[$\bullet$] $L^{-}(n)=\{A\in \FF^{n\times n} \,|\, AJ$ is a skew-symmetric
matrix$\}$.
\end{enumerate}

Denote coefficients in the characteristic polynomial
of an $n\times n$ matrix $X$ by $\si_t(X)$, i.e., %
$$\det(\la E-X)=\la^n-\si_1(X)\la^{n-1}+\cdots+(-1)^n\si_n(X).$$ %
So, $\si_1(X)=\tr(X)$ and $\si_n(X)=\det(X)$.

Assume $n$ is even. Define the {\it generalized pfaffian} of an arbitrary $n\times n$
matrix $X=(x_{ij})$ by
$$\P(X)=\Pf(X-X^T),$$
\noindent where $\Pf$ stands for the pfaffian of a skew-symmetric matrix. By abuse of
notation we will refer to $\P$ as the pfaffian. For $\FF=\QQ$ there is a more convenient
formula
\begin{eq}\label{eq_P}
\P(X)=\Pf(X-X^T)=\frac{1}{(n/2)!}\sum\limits_{\pi\in \Symmgr_{n}}\sign(\pi)
\prod\limits_{i=1}^{n/2} x_{\pi(2i-1),\pi(2i)}.%
\end{eq} %

For $n\times n$ matrices $X_1=(x_{ij}(1)),\ldots,X_s=(x_{ij}(s))$ and positive integers
$r_1,\ldots,r_s$, satisfying $r_1+\cdots+r_s=n/2$, consider the polynomial $\P(x_1
X_1+\cdots+x_s X_s)$ in the variables $x_l,\ldots,x_s$. The partial linearization %
$$\P_{r_1,\ldots,r_s}(X_1,\ldots,X_s)$$ %
of the pfaffian is the coefficient at $x_1^{r_1}\cdots x_s^{r_s}$ in this polynomial.
In other words, for $\FF=\QQ$ we have %
\begin{eq}\label{eq_lin_P}
{\P}_{r_1,\ldots,r_s}(X_1,\ldots,X_s)= \frac{1}{c}\sum\limits_{\pi\in
\Symmgr_{n}}\sign(\pi)
\prod\limits_{j=1}^s\;\prod\limits_{i=r_1+\cdots+r_{j-1}+1}^{r_1+\cdots+r_{j}}
x_{\pi(2i-1),\pi(2i)}(j),%
\end{eq} %
where $c=r_1!\cdots r_s!$. The partial linearization
$\det_{r_1,\ldots,r_s}(X_1,\ldots,X_s)$ of the determinant is defined analogously, where
$X_1,\ldots,X_s$ are $n\times n$ matrices, $r_1+\cdots+r_s=n$, and $n$ is arbitrary.

\subsection{Rings}\label{subsection_rings}
In what follows, ${\NN}$ stands for the set of non-negative integers.

Denote by
$\algA=\FF[x_1,\dots,x_d]$ %
the polynomial ring in $x_1,\ldots,x_d$ over $\FF$, i.e., $\algA$ is a commutative
$\FF$-algebra with unity generated by algebraically independent elements
$x_1,\ldots,x_d$.

%

Let $\{f_i\,|\,i\in I\}$ be a finite or countable system of generators for a commutative
algebra $\algA$. Then
$$\algA\simeq \FF[x_i\,|\,i\in I]/ T,$$
where $T$ is the ideal of {\it relations}.

\subsection{Quivers}\label{subsection_def_quivers}
A {\it quiver} $\Q=(\Q_0,\Q_1)$ is a finite oriented graph, where $\Q_0$ is the set of
vertices and $\Q_1$ is the set of arrows.  Multiple arrows and loops in $Q$ are allowed.
For an arrow $\al$, denote by $\al'$ its head
and by $\al''$ its tail, i.e., %
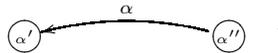
$$\vcenter{
\xymatrix@C=1cm@R=1cm{ %
\vtx{\al'}\ar@/^/@{<-}[rr]^{\al} && \vtx{\al''}\\
}} \quad.
$$
We say that $\al=\al_1\cdots \al_s$ is a {\it path} in $\Q$ (where $\al_1,\ldots,\al_s\in
\Q_1$), if $\al_1''=\al_2',\ldots,\al_{s-1}''=\al_s'$, i.e.,
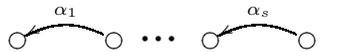
$$\vcenter{
\xymatrix@C=1cm@R=1cm{ %
\vtx{\;}\ar@/^/@{<-}[r]^{\al_1} & %
\vtx{\;} & %
\vtx{\;}\ar@/^/@{<-}[r]^{\al_s} & %
\vtx{\;}\\
}} \quad.
\begin{picture}(0,0)
\put(-73,-3){
\put(0,0){\circle*{2}} %
\put(-5,0){\circle*{2}} %
\put(5,0){\circle*{2}} %
} %
\end{picture}
$$
The head of the path $\al$ is $\al'=\al_1'$ and the tail is $\al''=\al_s''$.  A path
$\al$ is called {\it closed} if $\al'=\al''$. A closed path $\al$ is called {\it
incident} to a vertex $v\in\Q_0$ if $\al'=v$. Similarly, closed paths
$\be_1,\ldots,\be_s$ in $Q$ are called {\it incident} to $v$ if $\be_1'=\cdots=\be_s'=v$.
A closed path $\al$ is called {\it primitive} if $\al$ is not equal to a power of a shorter
closed path, i.e., $\al\neq \be^l$ for any $l>1$ and any closed path $\be$ in
$\Q$.

\section{Matrix invariants and representations of free algebras}\label{section_matrix}

\subsection{Generators for matrix invariants}\label{subsection_matrix}

Let $G$ be a group from the list $GL(n)$, $O(n)$, $\Sp(n)$, $SO(n)$ and let
$$H=\FF^{n\times n}\oplus\cdots\oplus \FF^{n\times n}$$
be $d$-tuple of $n\times n$ matrices over $\FF$. The group $G$ acts on $H$ by the
diagonal conjugation, i.e.,
\begin{eq}\label{eq_diag_conj}
g\cdot (A_1,\ldots,A_d)=(g A_1 g^{-1},\ldots,g A_d g^{-1}),
\end{eq}
where $g\in G$ and $A_1,\ldots,A_d\in \FF^{n\times n}$.

The coordinate ring of $H$ (i.e.~the ring of polynomial functions $f:H\to \FF$) is
the ring of polynomials  %
$$\FF[H]=\FF[x_{ij}(k)\,|\,1\leq i,j\leq n,\, 1\leq k\leq d],$$ %
where $x_{ij}(k)$ stands for the coordinate function on $H$ that takes a representation
$(A_1,\ldots,A_d)\in H$ to the $(i,j)^{\rm th}$ entry of the matrix $A_k$. Denote by
$$X_k=\left(\begin{array}{ccc}
x_{11}(k) & \cdots & x_{1n}(k)\\
\vdots & & \vdots \\
x_{n1}(k) & \cdots & x_{nn}(k)\\
\end{array}
\right)$$%
the $k^{\rm th}\!\!$ {\it generic} matrix ($1\leq k\leq n$).

The action of $G$ on $H$ induces the action on $\FF[H]$ as follows: $(g\cdot
f)(h)=f(g^{-1}\cdot h)$ for all $g\in G$, $f\in \FF[H]$,
$h\in H$. In other words,  %
$$g\cdot x_{ij}(k)= (i,j)^{\rm th}\text{ entry of }g^{-1}X_k g.$$%
The algebra of {\it matrix invariants} is
$$\FF[H]^{G}=\{f\in \FF[H]\,|\,g\cdot f=f\;{\rm for\; all}\;g\in G\}.$$
We do not consider the case $G=SL(n)$ because invariants for $GL(n)$ and $SL(n)$ are the
same.

\begin{theo}\label{theo_matrix}
The algebra of matrix invariants $\FF[H]^G$ is generated by the following elements:
\begin{enumerate}
\item[a)] $\si_t(A)$ ($1\leq t\leq n$ and $A$ ranges over all monomials in
$X_1,\ldots,X_d$), if $G=GL(n)$;

\item[b)] $\si_t(B)$ ($1\leq t\leq n$), if $G=O(n)$;

\item[c)] $\si_t(B)$ ($1\leq t\leq n$), if $G=SO(n)$ and $n$ is odd;

\item[d)] $\si_t(B)$, $\P_{r_1,\ldots,r_s}(B_1,\ldots,B_s)$ ($1\leq t\leq n$,
$r_1+\cdots+r_s=n/2$),  if $G=SO(n)$ and $n$ is even.
\end{enumerate}

In~b),~c), and~d) matrices $B,B_1,\ldots,B_s$ range over all monomials in
$X_1,\ldots,X_d$, $X_1^T,\ldots,X_d^T$.

\begin{enumerate}
\item[e)] $\si_t(C)$ ($1\leq t\leq n$ and $C$ ranges over all monomials in
$X_1,\ldots,X_d$, $JX_1^TJ,\ldots,JX_d^TJ$), if $G=\Sp(n)$.
\end{enumerate}
\end{theo}

\example\label{ex1} Let $G=GL(n)$, $n=2$. Then a
minimal (i.e.~irreducible) system of generators for the algebra $\FF[H]^G$ is the following one:
\begin{enumerate}
\item[$\bullet$] $\tr(X_{k_1}\cdots X_{k_p}),\;\;\det(X_k)$ ($1\leq p\leq 3$, $1\leq k_1<\cdots< k_p\leq d$ and $1\leq k\leq d$), if the characteristic of $\FF$ is not two; 

\item[$\bullet$] $\tr(X_{k_1}\cdots X_{k_p}),\;\;\det(X_k)$ ($1\leq k_1<\cdots< k_p\leq d$, and $1\leq k\leq d$), if the characteristic of $\FF$ is two.  
\end{enumerate}
In particular, the maximal degree of elements from any  generating system is equal to $3$ ($d$, respectively) in the first (the second, respectively) case. If we also have $d=2$, then in both cases $\FF[H]^G$ is the polynomial ring in %
$$\tr(X_1),\;\tr(X_2),\;\tr(X_1 X_2),\;\det(X_1),\;\det(X_2).$$

\remark\label{rem_tr} In the case of a characteristic zero field it is enough to take
traces instead of $\si_t$, $1\leq t\leq n$, in the formulation of
Theorem~\ref{theo_matrix}.

\remark\label{rem_new1} The system of generators from Theorem~\ref{theo_matrix} is
infinite but it can be shown that there exists a finite generating system. In this remark
we give references to some results concerning finite generating systems.

From now on we assume that $G=GL(n)$ and consider the algebra of matrix invariants. An
upper bound on degrees of elements of a minimal system of generators (m.s.g.) was given
in~\cite{Domokos_gen_2002} in terms of the nilpotency degree of a (unitary) relatively
free finitely generated algebra with the identity $x^n=0$.

For small $n$ the algebra of matrix invariants can be described more explicitly. A
m.s.g.~is known for $n=2$
(see~\cite{Sibirskii_1968},~\cite{Procesi_1984},~\cite{DKZ_2002}) and $n=3$
(see~\cite{Lopatin_Comm1},~\cite{Lopatin_Comm2}). In characteristic zero case a
m.s.g.~for $n=4$ and $d=2$ was described in~\cite{Drensky_Sadikova_4x4}.

By the Noether normalization the algebra of matrix invariants contains a {\it system of
parameters}, i.e., a set of algebraically independent elements generating a subalgebra
such that the algebra of matrix invariants is integral over this subalgebra. A system of
parameters for matrix invariants was constructed for $n=2$ and any $d$
(see~\cite{Teranishi_1986} and~\cite{DKZ_2002}), $n=3,4$ and $d=2$
(see~\cite{Teranishi_1986}), and $n=d=3$ (see~\cite{Lopatin_Sib}).

For more detailed introduction to finite generating systems for matrix invariants see
overviews~\cite{Formanek_1987} and~\cite{Formanek_1991} by Formanek. For recent
developments in characteristic zero see~\cite{Drensky_survey_2007} and in positive
characteristic see~\cite{DKZ_2002}.

\remark\label{rem_history1} For a field of characteristic zero generators for matrix
invariants of $G\in\{GL(n),O(n),\Sp(n)\}$ were described by Sibirskii
in~\cite{Sibirskii_1968} and Procesi in~\cite{Procesi76}. Developing ideas
from~\cite{Procesi76}, Aslaksen et al.~calculated generators for $G=SO(n)$
(see~\cite{Aslaksen95}).

The importance of characteristic-free approach to matrix invariants was pointed out by
Formanek in overview~\cite{Formanek_1991} (see also~\cite{Formanek_1987}). Relying on the
theory of modules with good filtrations (see~\cite{Donkin85}), Donkin described
generators for matrix $GL(n)$-invariants in~\cite{Donkin92a}. As regards the rest of
classical linear groups over a field of positive characteristic, the first results were
obtained by Zubkov. In~\cite{Zubkov99} he described generators for matrix $O(n)$- and
$\Sp(n)$-invariants. The proof is based on ideas from~\cite{Donkin92a} and a reduction to
invariants of $\QS$-mixed representations with $\g_v=GL$ for every vertex $v$ (see below
Section~\ref{section_mixed_quivers} for the definition). The reduction was performed by
means of Frobenius reciprocity. Generators for matrix $SO(n)$-invariants were described
by Lopatin in~\cite{Lop_so_inv}.

\remark{} Note that it is possible to define $O(n)$ and $SO(n)$ in characteristic two  
case. But in this case even generators for invariants of several vectors are not known
(for the latest developments see~\cite{Domokos_Frenkel}).

\subsection{Representations of free algebras}\label{subsection_alg}

Consider a free associative algebra %
$$\algA=\FF\LA x_1,\ldots,x_d\RA$$ %
freely generated by $x_1,\ldots,x_d$. Denote by $\rep{n}{\algA}$ the category of its
$n$-dimensional representations. For a representation $\varphi\in \rep{n}{\algA}$ we
assign the point $h_{\varphi}=(\varphi(x_1),\ldots,\varphi(x_d))$ from $H=(\FF^{n\times
n})^d$. Conversely, any point $h$ in $H$ determines some representation of $\algA$. Two
representations $\varphi$ and $\psi$ are isomorphic if and only if $h_{\varphi}$ becomes
$h_{\psi}$ under a basis change of $\FF^n$, i.e., $GL(n)\cdot h_{\varphi}=GL(n)\cdot
h_{\psi}$. Hence, there is one to one correspondence between classes of isomorphic
representations and orbits of $GL(n)$ on $H$.

\begin{prop}\label{prop_closed}
A representation $\varphi\in\rep{n}{\algA}$ is semisimple if and only if the orbit
$GL(n)\cdot h_{\varphi}$ is closed in $H$.
\end{prop}
\bigskip %

Suppose an algebraic group $G$ acts on the affine variety $H$.

\definition{of a categorical quotient} A pair $(Y, \pi_{Y})$, where $Y$ is an
algebraic variety and $\pi_{Y}$ is a morphism of $H$ into $Y$, is called a {\it
categorical quotient} for the action of $G$ on $H$ if for any morphism of algebraic
varieties $\xi:H \to Z$ that is constant on the orbits of $G$ there exists a unique
morphism $\eta: Y\to Z$ such that $\xi=\eta\circ \pi_Y$.
\bigskip

Denote by $H//G$ the set of all closed orbits of $G$ on $H$.

\begin{prop}\label{prop_quotient}
\begin{enumerate}
\item[a)] The set $H//G$ can be endowed with the structure of affine variety in such a
way that $H//G$ becomes a categorical quotient for the action of $G$ on $H$;

\item[b)] $\FF[H//G]=\FF[H]^{G}$. %
\end{enumerate}
\end{prop}

\begin{cor}\label{cor_quotient} Let $\varphi,\psi\in\rep{n}{\algA}$ be semisimple representations.
Then $\varphi\simeq\psi$ if and only if $f(h_{\varphi})=f(h_{\psi})$ for all
$f\in\FF[H]^G$.
\end{cor}

\subsection{Representations of free algebras with involution}\label{subsection_inv}
In this section we assume that $\FF$ is an {\it algebraically closed} field of
characteristic different from two. Then results from the previous section can be proven
for algebras with involution. Let us state it in a precise way.

An algebra $\algB$ with unity is called an {\it algebra with involution} $\ast$ if
\begin{enumerate}
\item[$\bullet$] $\ast:\algB\to \algB$ is a linear map;

\item[$\bullet$] $1^{\ast}=1$, where $1$ stands for the unity of $\algB$;

\item[$\bullet$] $(xy)^{\ast}=y^{\ast}x^{\ast}$ for all $x,y\in\algB$;

\item[$\bullet$] $x^{\ast\ast}=x$ for all $x\in\algB$.
\end{enumerate}
So we assume that $\ast$ is identical on $\FF\subset \algB$.

An $\varepsilon$-form (where $\varepsilon=\pm1$) on a vector space $V=\FF^n$ is a
bilinear nondegenerate form $\LA\cdot,\cdot\RA:V\otimes V\to\FF$ which is
$\varepsilon$-symmetric, i.e., $\LA u,v\RA=\varepsilon\LA v,u\RA$ for all $u,v\in V$.
Given an $\varepsilon$-form on $V$, the algebra $\End_{\FF}(V)$ can be endowed with a
{\it canonical} involution defined by the property: $\LA f^{\ast}(u),v\RA=\LA u,f(v)\RA$
for all $u,v\in V$ and $f\in \End_{\FF}(V)$. Since any nondegenerate symmetric
(skew-symmetric, respectively) bilinear form is isomorphic to the bilinear form defined
by the matrix $E$ ($J$, respectively), we can assume that a canonical involution on
$\End_{\FF}(V)$ is one of the following involutions:
\begin{enumerate}
\item[a)] $A\to A^T$,

\item[b)] $A\to -JA^TJ$,
\end{enumerate}
where $A\in \End_{\FF}(V)\simeq\FF^{n\times n}$. These involutions are called {\it
orthogonal} and {\it symplectic}, respectively.

An $n$-dimensional representation of the algebra with involution $\algB$ is a
$\varphi\in\rep{n}{\algB}$ (i.e., $\varphi:\algB\to \End_{\FF}(V)$, where $V=\FF^n$, is a
representation of $\algB$ considered as a usual algebra) such that
\begin{enumerate}
\item[$\bullet$] $\End_{\FF}(V)$ is endowed with {\it canonical} involution which we also
denote by $\ast$;

\item[$\bullet$] $\varphi(x^{\ast})=\varphi(x)^{\ast}$ for all $x\in \algB$.
\end{enumerate}
So we consider two kinds of representations: orthogonal and symplectic. Denote the
categories of these representations by $\repO{n}{\algB}$ and $\repSp{n}{\algB}$,
respectively. Let $G$ stands for $O(n)$ in the orthogonal case and for $\Sp(n)$ in the
symplectic case.

Consider a free associative algebra $\algB=\FF\LA x_1,\ldots,x_d,y_1,\ldots,y_d\RA$ and
endow it with involution as follows: $x_1^{\ast}=y_1,\ldots,x_d^{\ast}=y_d$. Then $\algB$
is called a free algebra with involution.

As in the previous section, for a $\varphi\in \repO{n}{\algB}$ we assign the point
$h_{\varphi}\in H=(\FF^{n\times n})^d$. Then we can see that there is one to one
correspondence between classes of isomorphic representations from $\repO{n}{\algB}$ and
orbits of $G$ on $H$. Moreover, analogues of
Propositions~\ref{prop_closed},~\ref{prop_quotient} and Corollary~\ref{cor_quotient} are
valid. Similar results are also valid for symplectic representations.

\subsection{The definition of $\si_{t,r}$}\label{subsection_DP}
Assume that $X,Y,Z$ are $n\times n$ matrices and $t,r\geq0$. To describe relations
between $O(n)$-matrix invariants, we need the polynomial $\si_{t,r}(X,Y,Z)$. In order to
define it, consider the quiver $\Q$
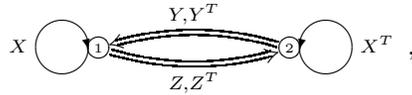
$$
\loopR{0}{0}{X} %
\xymatrix@C=1cm@R=1cm{ %
\vtx{1}\ar@2@/^/@{<-}[rr]^{Y,Y^T} &&\vtx{2}\ar@2@/^/@{<-}[ll]^{Z,Z^T}\\
}%
\loopL{0}{0}{X^T}\qquad,
$$
where there are two arrows from vertex $1$ to vertex $2$ and there are two arrows in the
opposite direction. By abuse of notation arrows of $\Q$ are denoted by the same letters
as matrices. Hence, any path in $\Q$ (see Section~\ref{subsection_def_quivers}) can be
interpreted as a product of matrices. Note that for any $\al\in \Q_1$ an arrow
$\al^T\in\Q_1$ is uniquely defined.

Let $\al=\al_1\cdots \al_p$ and $\be=\be_1\cdots \be_q$ be closed paths in $\Q$, where
$\al_1,\ldots,\al_p,\be_1,\ldots,\be_q\in \Q_1$. Then
\begin{enumerate}
\item[$\bullet$] we say that $\al$ and $\be$ are {\it equal} and write $\al=\be$ if $p=q$
and $\al_1=\be_1,\ldots,\al_p=\be_p$;

\item[$\bullet$] $\al^T$ is a closed path in $\Q$ defined by $\al^T=\al_p^T\cdots
\al_1^T$;

\item[$\bullet$] we say that $\al$ and $\be$ are {\it equivalent} and write $\al\sim\be$
if there exists a cyclic permutation $\pi\in\Symmgr_p$ such that
$\al_{\pi(1)}\cdots\al_{\pi(p)}=\be$ or $\al_{\pi(1)}\cdots\al_{\pi(p)}=\be^T$.
\end{enumerate}
As an example, we have $(X Y^T Z)^T=Z^T Y X^T$, $Y Z\sim Y^T Z^T$, $X Y X^T Z\sim X Y^T
X^T Z^T$.

Denote the degree of a path $\al$ in $\be\in\Q_1$ by $\deg_{\be}(\al)$ and the
multidegree of a path $\al$ by
$$\mdeg(\al)=%
(\deg_{X}(\al)+\deg_{X^T}(\al),\;\deg_{Y}(\al)+\deg_{Y^T}(\al),\;\deg_{Z}(\al)+\deg_{Z^T}(\al)).$$
As an example, $\deg_{Y}(Y Z)=1$, $\deg_{Y^T}(Y Z)=0$, and $\mdeg(X Y X^T X^T
Z^T)=(3,1,1)$.

\definition{of $\si_{t,r}(X,Y,Z)$} Let ${\mathcal P}$ be a set of
representatives of equivalence classes of primitive closed paths in $\Q$ (see Section~\ref{subsection_def_quivers} for the definition). Then
$$\si_{t,r}(X,Y,Z)=\sum %
(-1)^{t+\sum_{i=1}^p j_i(\deg_Y(\al_i)+\deg_Z(\al_i)+1)} \;%
\si_{j_1}(\al_1)\cdots\si_{j_p}(\al_p),$$ %
where the sum ranges over pairwise different $\al_1,\ldots,\al_p\in {\mathcal P}$ and
$j_1,\ldots,j_p\geq1$ such that
$$j_1\mdeg(\al_1)+\cdots+j_p\mdeg(\al_p)=(t,r,r).$$
For $t=r=0$ we define $\si_{0,0}(X,Y,Z)=1$.

\example\label{ex_DP} {\bf 1.} If $t=0$ and $r=1$, then %
$${\mathcal P}=\{Y Z,\, Y Z^T,\, \ldots\}.$$
Hence $\si_{0,1}(X,Y,Z)=-\tr(YZ)+\tr(YZ^T)$
\smallskip

{\bf 2.} If $t=r=1$, then %
$${\mathcal P}=\{X,\, Y Z,\, Y Z^T,\, X Y Z,\, X Y Z^T,\, X Y^T Z,\, X Y^T Z^T,\, \ldots\}$$
and we can see that $\si_{1,1}(X,Y,Z)=$
$$-\tr(X)\tr(YZ)+\tr(X)\tr(YZ^T)+ \\%
\tr(X Y Z)-\tr(X Y Z^T)-\tr(X Y^T Z)+\tr(X Y^T Z^T).$$ %

\remark\label{rem01} $\si_{t,0}(X,Y,Z)=\si_t(X)$.

\remark\label{rem_history001} The complete linearization of $\si_{t,r}$ was introduced by
Procesi (see~\cite{Procesi76}, Section~8 of Part~I). Then $\si_{t,r}$ was introduced
by Zubkov in~\cite{ZubkovII}. Note that we have $\si'_{t+2r,r}(X,Z,Y)=\si_{t,r}(X,Y,Z)$, where $\si'_{t+2r,r}$ stands for the function defined in~\cite{ZubkovII}.  The mentioned definitions from~\cite{Procesi76} and~\cite{ZubkovII} are different from our
definition, but their equivalence can be established applying the decomposition formula
from~\cite{Lop_bplp}.

Another way to define $\si_{t,r}$ is via the determinant-pfaffian $\DP_{r,r}(X,Y,Z)$
that was defined in~\cite{LZ1} as a ``mixture"{} of the determinant of $X$ and pfaffians
of $Y$ and $Z$. See below part~4 of Example~\ref{ex_pf_1} for the definition. One can
show that $\DP_{r,r}$ relates to $\si_{t,r}$ in the same way as the determinant relates
to $\si_t$, i.e., for  and $n\times n$ matrices $X,Y,Z$ we have
\begin{enumerate}
\item[$\bullet$] $\DP_{r,r}(X+\la E,Y,Z)=\sum_{t=0}^{t_0} \la^{t_0-t}\si_{t,r}(X,Y,Z)$, where $n=t_0+2r$, $t_0\geq0$;

\item[$\bullet$] $\det(X+\la E)=\sum_{t=0}^{n} \la^{n-t}\si_t(X)$.
\end{enumerate}
In particular, if $n=t+2r$, then $\DP_{r,r}(X,Y,Z)=\si_{t,r}(X,Y,Z)$. Note that this
approach gives us $\si_{t,r}(X,Y,Z)$ as a polynomial in entries of matrices $X,Y,Z$. But
for our purposes we have to present $\si_{t,r}(X,Y,Z)$ in a different way, namely, as a
polynomial in $\si_t(\alpha)$, where $t$ ranges over positive integers and $\al$ ranges
over monomials in $X,Y,Z,X^T,Y^T,Z^T$.

\subsection{Relations for matrix invariants}\label{subsection_relations}

In this section we use notations from Section~\ref{subsection_matrices}. We assume that
$G$ is $GL(n)$ or $O(n)$. We denote by $\M$ the monoid freely generated by letters  $x_1,\ldots,x_d$ ($x_1,\ldots,x_d,x_1^T,\ldots,x_d^T$, respectively) if $G=GL(n)$
($G=O(n)$, respectively). Given $\al=\al_1\cdots \al_s\in\M$, where $\al_i=x_i$ or $\al_i=x_i^T$, we assume $X_{\al}=X_{\al_1}\cdots X_{\al_s}$, where 
$$X_{\al_i}=\left\{
\begin{array}{ccl}
X_i  &,&\text{if }\al_i=x_i\\
X_i^T&,&\text{if }\al_i=x_i^T\\
\end{array}
\right..
$$ %
Let us recall that $\al\in\M$ is {\it primitive} if $\al$ is not equal to a power of a shorter monomial.

\begin{lemma}\label{lemma_primitive}
\begin{enumerate}
\item[a)] The algebra $\FF[H]^G$ is generated by $\si_t(X_{\al})$, where $1\leq t\leq n$ and
$\al\in\M$ is primitive.

\item[b)] For $a_i\in\FF$, $\al_i\in\M$, and $t>0$ the element $\si_t(\sum_i a_iX_{\al_i})$ is a polynomial in $\si_l(X_{\al})$, where $l$ ranges over $\{1,\ldots,n\}$ and $\al$ ranges over primitive elements from $\M$.

\item[c)] Let $G=O(n)$, $a_i,b_j,c_k\in\FF$, $\al_i,\be_j,\ga_k\in\M$, and $t,r\geq0$. Then the
element $\si_{t,r}(\sum_i a_i X_{\al_i},\,\sum_j b_j X_{\be_j},\,\sum_k c_k X_{\ga_k})$ is a polynomial in $\si_l(X_{\al})$, where $l$ ranges over $\{1,\ldots,n\}$ and $\al$ ranges over primitive elements from $\M$.
\end{enumerate}
\end{lemma}

\begin{proof}
In this proof $A,A_1,\ldots,A_p$ stands for an arbitrary $n\times n$ matrices.

a) For $1\leq t\leq n$ and $l\geq2$ we have the following well-known formula:
\begin{eq}\label{eq_D}
\si_t(A^l)=\sum\limits_{i_1,\ldots,i_{t l}\geq0}b^{(t,l)}_{i_1,\ldots,i_{t l}}
    \si_1(A)^{i_1}\cdots\si_{t l}(A)^{i_{t l}},%
\end{eq}
\noindent where we assume that $\si_i(A)=0$ for $i>n$. 
In~\Ref{eq_D} coefficients $b^{(t,l)}_{i_1,\ldots,i_{rl}} \in \ZZ$ do not depend on $A$
and $n$. If we take $A=\diag(a_1,\ldots,a_n)$ is a diagonal matrix, then $\si_t(A^l)$ is
a symmetric polynomial in $a_1,\ldots,a_n$ and $\si_i(A)$ for $1\leq i\leq n$ is the
$i^{\rm th}$ elementary symmetric polynomial in $a_1,\ldots,a_n$. Thus, the coefficients
$b^{(t,l)}_{i_1,\ldots,i_{rl}}$ with $rl\leq n$ can easily be found. As an example,
\begin{eq}\label{eq_tr_a2}
\tr(A^2)=\tr(A)^2-2\si_2(A).
\end{eq}
Using~\Ref{eq_D} we complete the proof. %
\medskip

b) We start with a definition. Let $y_1,\ldots,y_p$ be some letters. Then a {\it cycle}
in these letters is an equivalent class of some monomial in $y_1,\ldots,y_p$ with respect
to cyclic permutations.

For $1\leq t\leq n$ Amitsur's formula states~\cite{Amitsur_1980}:
\begin{eq}\label{eq_Amitsur}
\si_t(A_1+\cdots+A_p)=\sum (-1)^{t-(j_1+\cdots+j_q)} \si_{j_1}(c_1)\cdots\si_{j_q}(c_q),
\end{eq}
where the sum ranges over all pairwise different primitive cycles $c_1,\ldots,c_q$ in
letters $A_1,\ldots,A_p$ and positive integers $j_1,\ldots,j_q$ with
$\sum_{i=1}^{q}j_i\deg(c_i)=t$. As an example,
$$
\si_2(A_1+A_2)=\si_2(A_1)+\si_2(A_2)+\tr(A_1)\tr(A_2)-\tr(A_1A_2).
$$%
Using~\Ref{eq_D} together with the formula
\begin{eq}
\si_t(a A)=a^t\si_t(A),
\end{eq}
where $a\in\FF$, we complete the proof. %
\medskip

c) Follows from part~b).
\end{proof}

For $a_i,b_j,c_k\in\FF$, $\,\al,\be,\al_i,\be_j,\ga_k\in\M$, and $t,r\geq0$ we consider the
following equalities:
\begin{enumerate}
\item[a)] $\si_t(X_{\al}X_{\be})-\si_t(X_{\be}X_{\al})=0$ for $1\leq t\leq n$;

\item[b)] $\si_t(\sum_i a_iX_{\al_i})=0$ for $t>n$;

\item[c)] $\si_{t,r}(\sum_i a_i X_{\al_i},\,\sum_j b_j X_{\be_j},\,\sum_k c_k X_{\ga_k})=0$ for $t+2r>n$.

\item[d)] free relations, i.e., relations between $O(m)$-invariants that are valid for any $m>0$.
\end{enumerate}
Here we assume that Lemma~\ref{lemma_primitive} is applied to these equalities, i.e., the
left hand sides of the above equalities are polynomials in $\si_l(X_{\al})$, where $l$ ranges over $\{1,\ldots,n\}$ and $\al$ ranges over primitive elements from $\M$.

\begin{theo}\label{theo_relations} Consider the system of generators for $\FF[H]^G$
from part~a) of Lemma~\ref{lemma_primitive}. Then the ideal of relations between them is
generated by the following relations:
\begin{enumerate}
\item[$\bullet$] a) and b), if $G=GL(n)$;

\item[$\bullet$] a), c), and d), if $G=O(n)$.
\end{enumerate}
\end{theo}

\remark\label{rem2} Remark~\ref{rem01} implies that relation c) with $r=0$ is the same as
relation b).

\remark\label{rem_history2} In case of characteristic zero relations for matrix $GL(n)$-,
$O(n)$- and $\Sp(n)$-invariants were described by Procesi in~\cite{Procesi76}.
Independently, relations for $G=GL(n)$ were established by Razmyslov 
in~\cite{Razmyslov74}. Over a field of arbitrary characteristic relations for matrix
$GL(n)$-invariants were described by Zubkov (see~\cite{Zubkov96}) and relations for
matrix $O(n)$-invariants were described by Lopatin
(see~\cite{Lop_relations}). The last result is based on the reduction to mixed
representations of quivers that was performed by Zubkov in~\cite{ZubkovII}.


\section{Representations of quivers}\label{section_quivers} %
\subsection{Definitions}\label{subsection_quivers_definitions} %

The notion of a representation of a quiver generalizes the notion of several linear maps
on a vector space (see Section~\ref{subsection_matrix}).

Consider a quiver $\Q=(\Q_0,\Q_1)$ and assume that $\Q_0=\{1,\ldots,l\}$ for some $l$.
Given a {\it dimension vector} $\n=(\n_1,\ldots,\n_l)$, we assign an $\n_v$-dimensional
vector space $V_v$ to $v\in \Q_0$. We identify $V_v$ with the space of column vectors
$\FF^{\n_v}$. Fix the {\it standard} basis $e(v,1),\ldots,e(v,\n_v)$ for $\FF^{\n_v}$,
where $e(v,i)$ is a column vector whose $i^{\rm th}$ entry is $1$ and the rest of entries
are zero. A {\it representation} of $\Q$ of dimension vector
$\n$ is a collection of matrices %
$$h=(h_\al)_{\al\in \Q_1}\in %
H=H(\Q,\n)=\bigoplus_{\al\in \Q_1} \FF^{\n_{\al'}\times \n_{\al''}} \simeq%
\bigoplus_{\al\in \Q_1} \Hom_\FF(V_{\al''},V_{\al'}),$$ %
where $\FF^{n_1\times n_2}$ stands for the linear space of $n_1\times n_2$ matrices over
$\FF$ and the isomorphism is given by the choice of bases.
The action of the group %
$$G=GL(\n)=\prod_{v\in \Q_0} GL(\n_v)$$ %
on $H$ is via the change of the bases for $V_v$ ($v\in \Q_0$). In
other words, $GL(\n_v)$ acts on $V_v$ by left multiplication, and
this action induces the action of $G$ on $H$ by
$$g\cdot h=(g_{\al'}h_\al g_{\al''}^{-1})_{\al\in \Q_1},$$ %
where $g=(g_\al)_{\al\in \Q_1}\in G$ and $h=(h_\al)_{\al\in \Q_1}\in H$. We refer to the
pair $(\Q,\n)$ as the {\it quiver setting}.

The coordinate ring of $H$ (i.e.~the ring of polynomial functions $f:H\to\FF$)
is the polynomial ring  %
$$\FF[H]=\FF[x_{ij}(\al)\,|\,\al\in \Q_1,\,1\leq i\leq \n_{\al'},1\leq j\leq\n_{\al''}].$$ %
Here $x_{ij}(\al)$ stands for the coordinate function on $H$ that takes $h\in H$ to the
$(i,j)^{\rm th}$ entry of the matrix $h_{\al}$. Denote by %
$$X_\al=\left(\begin{array}{ccc}
x_{1,1}(\al) & \cdots & x_{1,\n_{\al''}}(\al)\\
\vdots & & \vdots \\
x_{\n_{\al'},1}(\al) & \cdots & x_{\n_{\al'},\n_{\al''}}(\al)\\
\end{array}
\right)$$%
the $\n_{\al'}\times \n_{\al''}$ {\it generic} matrix.

The action of $G$ on $H$ induces the action on $\FF[H]$ as follows: $(g\cdot
f)(h)=f(g^{-1}\cdot h)$ for all $g\in G$, $f\in \FF[H]$,
$h\in H$. In other words, %
$$g\cdot x_{ij}(\al)=(i,j)^{\rm th}\text{ entry of }g_{\al'}^{-1}X_\al g_{\al''}$$ %
for all $g\in G$, $\al\in \Q_1$. The algebra of {\it invariants} is
$$\FF[H]^{G}=\{f\in \FF[H]\,|\,g\cdot f=f\;{\rm for\; all}\;g\in G\}.$$

\example\label{ex01} Let $\Q$ be a quiver
$$\vcenter
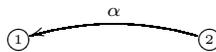
{
\xymatrix@C=1cm@R=1cm{ %
\vtx{1}\ar@/^/@{<-}[rr]^{\al} &&\vtx{2}\\
}} \quad
$$
and $\n=(n,m)$. Then the group $G=GL(n)\times GL(m)$ acts on $H=\FF^{n\times m}$ by the
rule
$$(g_1,g_2)\cdot A=g_1 A g_2^{-1},$$
where $(g_1,g_2)\in G$ and $A\in \FF^{n\times m}$. Hence the orbits of this action
correspond to linear maps from $\FF^{m}$ to $\FF^{n}$.

\example\label{ex02} Let $\Q$ be a quiver with one vertex and $d$ loops:
$$
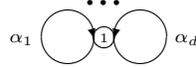
\loopR{0}{0}{\al_1}%
\xymatrix@C=1cm@R=1cm{ %
\vtx{1}\\
}\loopL{0}{0}{\al_d}
\begin{picture}(0,20)
\put(-25,16){\circle*{2}} %
\put(-20,16){\circle*{2}} %
\put(-30,16){\circle*{2}} %
\end{picture}
$$
and $\n=(n)$. Then the group $G=GL(n)$ acts on $H=(\FF^{n\times n})^d$ by the diagonal
conjugation (see~\Ref{eq_diag_conj}).  Obviously, the orbits of this action correspond to
$d$-tuples of linear maps on $\FF^{n}$. Generators for $\FF[H]^G$ were described in
part~a) of Theorem~\ref{theo_matrix}.

\example\label{ex03} Let $\Q$ be a quiver
$$\vcenter{
\xymatrix@C=1cm@R=1cm{ %
\vtx{1}\ar@/^/@{<-}[r]^{\al} &\vtx{2}\ar@/^/@{<-}[r]^{\be}&\vtx{3}\\
}} \quad
$$
and $\n=(n,m,l)$. Then the group $G=GL(n)\times GL(m)\times GL(l)$ acts on
$H=\FF^{n\times m}\oplus \FF^{m\times l}$ by the rule
$$(g_1,g_2,g_3)\cdot (A,B)=(g_1 A g_2^{-1}, g_2 B g_3^{-1}).$$

\subsection{Generators for invariants of representations of quivers} \label{subsection_quivers_generators}
Let us recall that the notion of closed path was given in
Section~\ref{subsection_def_quivers}.

\begin{theo}\label{theo_Donkin} Let $(\Q,\n)$ be a quiver setting. Then $\FF$-algebra
$\FF[H(\Q,\n)]^{GL(\n)}$ is generated by %
$$\si_t(X_{\be_1}\cdots X_{\be_s}),$$ %
where $\be_1\cdots \be_s$ is a closed path in $\Q$ and $1\leq t\leq \n_{\be'_1}$.
\end{theo}
\bigskip %

This theorem implies that if there is no closed paths in $\Q$, then
$$\FF[H(\Q,\n)]^{GL(\n)}=\FF.$$

\remark\label{rem_new2} The system of generators from Theorem~\ref{theo_Donkin} is
infinite but it can be shown that there exists a finite generating system. As an example,
invariants of quivers of dimension $(2,\ldots,2)$ were considered in~\cite{Lopatin_2222},
where an upper bound on degrees of elements of a minimal system of generators was given
and its precision was estimated.

\remark\label{rem_history3} Over a field of characteristic zero,
Theorem~\ref{theo_Donkin} was proven by Le Bruyn and Procesi
in~\cite{Le_Bruyn_Procesi_90}; over a field of arbitrary characteristic, it was proven by
Donkin~\cite{Donkin94}.

\subsection{Relations for invariants of representations of quivers} \label{subsection_quivers_relations}
Relations for invariants of quivers can be described in a similar way as relations for
matrix $GL(n)$-invariants (see Section~\ref{subsection_relations}). In what follows we will need the notions of primitive path and incident path that were introduced in Section~\ref{subsection_def_quivers}.

We assume that $(\Q,\n)$ is a quiver setting. Let $\M$ be the set of closed paths in $\Q$. For $\be=\be_1\cdots \be_s\in\M$, denote $X_{\be}=X_{\be_1}\cdots X_{\be_s}$. The following lemma is an analogue of Lemma~\ref{lemma_primitive}. 

\begin{lemma}\label{lemma_primitive2}
\begin{enumerate}
\item[a)] The algebra $\FF[H(\Q,\n)]^{GL(\n)}$ is generated by $\si_t(X_{\be})$, where   $\be=\be_1\cdots \be_s\in\M$ is primitive and $1\leq t\leq \n_{\be'_1}$.

\item[b)] Let $a_i\in\FF$, $\al_i\in\M$ be incident to $v\in\Q_0$, and $t>0$. Then the
element $\si_t(\sum_i a_iX_{\al_i})$ is a polynomial in $\si_l(X_{\al})$, where $l$ ranges over $\{1,\ldots,\n_v\}$ and $\al$ ranges over primitive elements from $\M$.
\end{enumerate}
\end{lemma}
\bigskip %

Given a vertex $v\in\Q_0$, we consider the following equalities:
\begin{enumerate}
\item[a)] $\si_t(X_{\al}X_{\be})-\si_t(X_{\be}X_{\al})=0$ for $1\leq t\leq \n_v$, where $\al,\be$ are paths in $\Q$ such that $\al\be\in\M$ is incident to $v$;

\item[b)] $\si_t(\sum_i a_iX_{\al_i})=0$ for $t>\n_v$, where $a_i\in\FF$ and 
$\al_i\in\M$ is incident to $v$.
\end{enumerate}
Here we assume that Lemma~\ref{lemma_primitive2} is applied to these equalities, i.e.,
the left hand sides of the above equalities are polynomials in $\si_l(X_{\al})$, where $l$
ranges over $\{1,\ldots,\n_v\}$ and $\al$ ranges over primitive elements from $\M$.

\begin{theo}\label{theo_quivers_relations} Consider the system of generators for $\FF[H(\Q,\n)]^{GL(\n)}$
from part~a) of Lemma~\ref{lemma_primitive2}. Then the ideal of relations between them is
generated by relations a) and b) considered for all vertices of $\Q$.
\end{theo}

\remark{} Theorem~\ref{theo_quivers_relations} was proved by Zubkov
in~\cite{Zubkov_Fund_Math_01}.

\section{$\QS$-mixed representations of quiver}\label{section_mixed_quivers}

\subsection{Definitions}\label{subsection_mixed_quivers_definitions}
The notion of representations of quivers can be generalized by the successive realization
of the following steps. At the end of this procedure we obtain {\it mixed}
representations of quivers.
\begin{enumerate}
\item[1.] Instead of $GL(\n)$ we can take a product $G(\n,\g)$ of classical linear
groups. Here $\g=(\g_1,\ldots,\g_l)$ is a vector, whose entries $\g_1,\ldots,\g_l$ are
symbols from the list $GL,O,\Sp,SL,SO$. By definition,
$$G(\n,\g)=\prod_{v\in \Q_0} G_v,$$
where
$$
G_v=\left\{%
\begin{array}{ll}
GL(\n_v),& \text{if }\g_v=GL\\
O(\n_v),& \text{if }\g_v=O\\
\Sp(\n_v),& \text{if }\g_v=\Sp\\
SL(\n_v),& \text{if }\g_v=SL\\
SO(\n_v),& \text{if }\g_v=SO\\
\end{array}
\right..
$$

Obviously, we have to assume that $\n$ and $\g$ are subject to the following
restrictions:
\begin{enumerate}
\item[a)] if $\g_v=\Sp$ ($v\in \Q_0$), then $\n_v$ is even;

\item[b)] if $\g_v$ is $O$ or $SO$ ($v\in \Q_0$), then the characteristic of $\FF$ is not
$2$.
\end{enumerate} %

\item[2.] We can change the definition of $G(\n,\g)$ in such a manner that allows us to
deal with bilinear forms together with linear mappings. Since bilinear forms on some
vector space $V$ are in one to one correspondence with linear mappings from the dual
vector space $V^{\ast}$ to $V$, we should change vector spaces assigned to some vertices
to the dual ones. In order to do this consider a mapping $\i:\Q_0\to \Q_0$ such that
\begin{enumerate}
\item[c)] $\i$ is an involution, i.e., $\i^2$ is the identical mapping;

\item[d)] $\n_{\i(v)}=\n_{v}$ for every vertex $v\in \Q_0$.
\end{enumerate}
For every $v\in \Q_0$ with $v<\i(v)$ assume that $V_{\i(v)}=V_v^{\ast}$. Consider the
dual basis $e(v,1)^{\ast},\ldots,e(v,\n_v)^{\ast}$ for $V_v^{\ast}$ and identify
$V_v^{\ast}$ with the space of column vectors of length $\n_v$, so $e(v,i)^{\ast}$ is the
same column vector as $e(v,i)$.

The action of $GL(\n_v)$ on $V_v$ induces the action on $V_v^{\ast}$, which we consider
as the degree one homogeneous component of the graded algebra $\FF[V_v]$. Given $g_v\in
GL(\n_v)$
and $u\in V_v^{\ast}$, we have %
$$g_v\cdot u=(g_v^{-1})^T u.$$
Hence, we should change the group $G$ to
$$G(\n,\g,\i)=\{g\in G(\n,\g)\,|\,
g_{\i(v)}=(g_v^{-1})^T\; {\rm for\;all}\; v\in \Q_0\; {\rm with}\;v<\i(v)\}.
$$
Since the vector spaces $\FF^n$ and $(\FF^{n})^{\ast}$ are isomorphic as modules over
$O(n)$, $\Sp(n)$, and $SO(n)$, we assume that
\begin{enumerate}
\item[e)] if $\g_v$ is $O$, $\Sp$ or $SO$ ($v\in \Q_0$), then $\i(v)=v$.
\end{enumerate}

\item[3.] Instead of the space $H(\Q,\n)$ we should take its subspace $H(\Q,\n,\h)$,
where $\h=(\h_\al)_{\al\in \Q_1}$ and $\h_\al$ is a symbol from the list
$M,S^{+},S^{-},L^{+},L^{-}$. By definition,  %
$$H(\Q,\n,\h)=\bigoplus_{\al\in \Q_1} H_{\al},$$
where
$$
H_{\al}=\left\{%
\begin{array}{ll}
\FF^{\n_{\al'}\times\n_{\al''}},& \text{if }\h_{\al}=M\\
S^{+}(\n_{\al'}),& \text{if }\h_{\al}=S^{+}\\
S^{-}(\n_{\al'}),& \text{if }\h_{\al}=S^{-}\\
L^{+}(\n_{\al'}),& \text{if }\h_{\al}=L^{+}\\
L^{-}(\n_{\al'}),& \text{if }\h_{\al}=L^{-}\\
\end{array}
\right..
$$
Additionally, we have to assume that $\n$ and $\h$ are subject to the
restriction: %
\begin{enumerate}
\item[f)] if $\h_\al\neq M$ ($\al\in \Q_1$), then $\n_{\al'}=\n_{\al''}$.
\end{enumerate}
\end{enumerate}
Consider a group $G=G(\n,\g,\i)\subset GL(\n)$ and a vector space $H=H(\Q,\n,\h)\subset
H(\Q,\n)$ satisfying the previous conditions~a)--f). To ensure that these inclusions
induce the action of $G$ on $H$, we assume that the following additional conditions are
also valid for all $v\in \Q_0$, $\al\in \Q_1$:
\begin{enumerate}\item[]
\begin{enumerate}
\item[g)] if $\al$ is a loop, i.e., $\al'=\al''$, and $\h_\al$ is $S^{+}$ or $S^{-}$,
then $\g_{\al'}$ is $O$ or $SO$;

\item[h)] if $\al$ is a loop and $\h_\al$ is $L^{+}$ or $L^{-}$, then $\g_{\al'}=\Sp$;

\item[i)] if $\al$ is not a loop and $\h_\al\neq M$, then $\i(\al')=\al''$ and $\h_\al$
is $S^{+}$ or $S^{-}$.
\end{enumerate}
\end{enumerate}
A quintuple $\QS=(\Q,\n,\g,\h,\i)$ satisfying~a)--i) is called a {\it mixed quiver
setting} and elements of $H$ are called {\it $\QS$-mixed representations} of the quiver $\Q$.
Definitions of the generic matrices $X_\al$ and the algebra of invariants $\FF[H]^G$ are
the same as above. Note
that %
if $\h_\al=S^{+}$, then $X_{\al}^T=X_\al$; %
if $\h_\al=S^{-}$, then $X_{\al}^T=-X_\al$; %
if $\h_\al=L^{+}$, then $(X_{\al}J)^T=X_\al J$; and %
if $\h_\al=L^{-}$, then $(X_{\al}J)^T=-X_\al J$. %

\example\label{ex3} \textbf{ 1.} Let $\Q$ be the following quiver
$$\vcenter
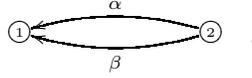
{
\xymatrix@C=1cm@R=1cm{ %
\vtx{1}\ar@/^/@{<-}[rr]^{\al} \ar@/_/@{<-}[rr]_{\be}&&\vtx{2}\\
}} \quad.
$$
Define a mixed quiver setting $\QS=(\Q,\n,\g,\h,\i)$ by $\n_1=\n_2=n$, $\g_1=\g_2=GL$,
$\h_\al=\h_\be=S^{+}$, and $\i(1)=2$. The group $G(\n,\g,\i)\simeq GL(n)$ acts on
$H(\Q,\n,\h)=S^{+}(n)\oplus S^{+}(n)$ by the rule
$$g\cdot (A,B)=(g A g^T, g B g^T)$$
for $g\in GL(n)$ and $(A,B)\in S^{+}(n)\oplus S^{+}(n)$. Hence the orbits of this action
correspond to pairs of symmetric bilinear forms on $\FF^{n}$. If we put
$\h_\al=\h_\be=S^{-}$, then we obtain pairs of skew-symmetric bilinear forms on
$\FF^{n}$. The classification problem for such pairs is a classical topic going back to
Weierstrass and Kronecker (see~\cite{Hodge_Pedoe_I} and~\cite{Gantmacher}).
\smallskip

\textbf{ 2.} Let $\Q$ be the following quiver
$$
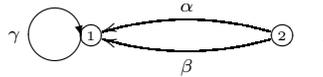
\loopR{0}{0}{\ga}%
\xymatrix@C=1cm@R=1cm{ %
\vtx{1}\ar@/^/@{<-}[rr]^{\al} \ar@/_/@{<-}[rr]_{\be}&&\vtx{2}\\
} \quad.
$$
Define $\QS=(\Q,\n,\g,\h,\i)$ by $\n_1=\n_2=n$, $\g_1=\g_2=GL$, $\h_\al=\h_\be=\h_\ga=M$,
and $\i(1)=2$. The group $G(\n,\g,\i)\simeq GL(n)$ acts on $H(\Q,\n,\h)=(\FF^{n\times
n})^3$ by the rule
$$g\cdot (A,B,C)=(g A g^T, g B g^T, g C g^{-1}).$$
Hence the orbits of this action correspond to pairs of bilinear forms on $\FF^{n}$
together with a linear map on $\FF^{n}$.
\smallskip

\textbf{ 3.} Let $\Q$ be the following quiver %
$$\vcenter
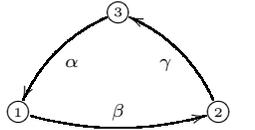
{
\xymatrix@C=1cm@R=1cm{ %
&\vtx{3}\ar@/_/@{->}[dl]^{\al} \ar@/^/@{<-}[dr]_{\ga}&\\
\vtx{1}\ar@/_/@{->}[rr]^{\be} &&\vtx{2}\\
}} \quad.
$$
Define a mixed quiver setting $\QS=(\Q,\n,\g,\h,\i)$ by $\n_1=\n_2=n$, $\n_3=m$;
$\g_1=\g_2=GL$, $\g_3=O$; $\h_\al=\h_\ga=M$, $\h_\be=S^{+}$; and $\i(1)=2$, $\i(3)=3$.
Hence the action of $G(\n,\g,\i)\simeq GL(n)\times O(m)$ on $H(\Q,\n,\h)=\FF^{n\times
m}\oplus S^{+}(n)\oplus \FF^{m\times n}$ is given by
$$(g_1,g_3)\cdot (A,B,C)=(g_1 A g_3^T, (g_1^{-1})^T B g_1^{-1}, g_3 C g_1^T)$$
for $(g_1,g_3)\in GL(n)\times O(m)$ and $(A,B,C)\in \FF^{n\times m}\oplus S^{+}(n)\oplus
\FF^{m\times n}$. %
\bigskip

Let us consider some partial cases of the notion of $\QS$-mixed representations. 

\example\label{ex3part4} We assume that $\QS=(\Q,\n,\g,\h,\i)$ is a mixed quiver setting.
\begin{enumerate}
\item[$\bullet$] If $\g_v=GL$ and $\i(v)=v$ for all $v\in\Q_0$, then $\h_{\al}=M$ for all
$\al\in\Q_1$; hence in this case $\QS$-mixed representations are the usual
representations of $\Q$.

\item[$\bullet$] If $\g_v\in\{GL,O\}$ for all $v\in\Q_0$ and $\h_{\al}\in\{M,S^{-}\}$ for
all $\al\in\Q_1$, then $\QS$-mixed representations are {\it orthogonal} representations
from~\cite{DW02}.

\item[$\bullet$] If $\g_v\in\{GL,\Sp\}$ for all $v\in\Q_0$ and $\h_{\al}\in\{M,S^{+}\}$
for all $\al\in\Q_1$, then $\QS$-mixed representations are {\it symplectic}
representations from~\cite{DW02}.

\item[$\bullet$] If $\g_v=GL$ for all $v\in\Q_0$ and $\h_{\al}=M$ for all $\al\in\Q_1$,
then $\QS$-mixed representations are {\it mixed} representations
from~\cite{Zubkov_preprint},~\cite{ZubkovI}.

\item[$\bullet$] If $\g_v\in\{GL,O,\Sp\}$ for all $v\in\Q_0$, then $\QS$-mixed representations are
{\it supermixed} representations of quivers from~\cite{Zubkov_preprint},~\cite{ZubkovI}, or
equivalently, representations of {\it signed} quivers from~\cite{Shmelkin}.
\end{enumerate}

\example\label{ex_semi} Let $\QS=(\Q,\n,\g,\h,\i)$ be a mixed quiver setting. We put  $G=G(\n,\g,\i)$ and $H=H(\Q,\n,\h)$.
\begin{enumerate}
\item[$\bullet$] If $\g_v=SL$ and $\i(v)=v$ for all $v\in\Q_0$, then $\h_{\al}=M$ for all
$\al\in\Q_1$. In this case $\FF[H]^G$ is the algebra of {\it semi-invariants} of  (usual) representations of $\Q$.

\item[$\bullet$] If $\g_v=SL$ for all $v\in\Q_0$ and $\h_{\al}=M$ for all $\al\in\Q_1$,
then $\FF[H]^G$ is the algebra of {\it semi-invariants} of mixed representations of $\Q$.

\item[$\bullet$] If $\g_v\in\{SL,SO\}$ for all $v\in\Q_0$ and
$\h_{\al}\in\{M,S^{+},S^{-}\}$ for all $\al\in\Q_1$, then $\FF[H]^G$ is the algebra of {\it semi-invariants} of supermixed representations of $\Q$.
\end{enumerate}

\subsection{Generators for invariants of supermixed representations of quivers}
\label{subsection_mixed_quivers_generators}

At first we consider an arbitrary mixed quiver setting and present some definitions. Then we restrict the general case to supermixed representations of quivers (see Theorem~\ref{theo_ZubI}). The general case is considered in  Section~\ref{subsection_general_case}.

Let $\QS=(\Q,\n,\g,\h,\i)$ be a mixed quiver setting and $\Q_0=\{1,\ldots,l\}$. Without
loss of generality we can assume that
\begin{eq}\label{eq_condition}
{\rm\; if\;} v\in \Q_0 {\rm\; and\;} \g_v=GL\text{ or }SL, {\rm\; then\;} \i(v)\neq v.
\end{eq}%
Otherwise we can add a new vertex $\ov{v}$ to $\Q$, and set $\i(v)=\ov{v}$,
$\n_{\ov{v}}=\n_v$, $\g_{\ov{v}}=\g_{v}$; this construction changes neither the space
$H(\Q,\n,\h)$ nor the algebra of invariants.

\definition{of the mixed double quiver setting $\QS^{\D}$}\label{def_double} %
Define the mixed {\it double} quiver setting $\QS^{\D}=(\Q^{\D},\n,\g,\h^{\D},\i)$
as follows: %
$$\Q_0^{\D}=\Q_0,\;\; \Q_1^{\D}=\Q_1\coprod \{\al^T\,|\,\al\in \Q_1,\;
\h_\al=M\}, $$ %
where $(\al^T)'=\i({\al''})$, $(\al^T)''=\i({\al'})$, and
$\h_{\al^T}^{\D}=M$ for $\al\in \Q_1$ with $\h_\al=M$ and $\h_{\al}^{\D}=\h_{\al}$ for
all $\al\in\Q_1$.

Define a mapping $\Phi^{\D}:\FF[H(\Q^{\D},\n,\h^{\D})]\to \FF[H(\Q,\n,\h)]$ such that
\begin{enumerate}
\item[a)] for $\al\in \Q_1$ we have $\Phi^{\D}(X_{\al})=X_\al$;

\item[b)] for $\al\in \Q_1$ and $\h_\al=M$ we define $\Phi^{\D}(X_{\al^T})$ as follows:
\begin{enumerate} %
\item[$\bullet$] If $\g_{\al'}\neq \Sp$ and $\g_{\al''}\neq \Sp$, then
$\Phi^{\D}(X_{\al^T})=X_\al^T$.

\item[$\bullet$] If $\g_{\al'}=\Sp$ and $\g_{\al''}\neq \Sp$, then
$\Phi^{\D}(X_{\al^T})=X_\al^T J(\n_{\al'})$.

\item[$\bullet$] If $\g_{\al'}\neq \Sp$ and $\g_{\al''}=\Sp$, then
$\Phi^{\D}(X_{\al^T})=J(\n_{\al''})X_\al^T$.

\item[$\bullet$] If $\g_{\al'}=\Sp$ and
$\g_{\al''}=\Sp$, then %
$\Phi^{\D}(X_{\al^T})=J(\n_{\al''}) X_\al^T J(\n_{\al'}) $.
\end{enumerate} %
\end{enumerate} %
Here $\Phi^{\D}(X_\al)$ stands for the matrix, whose $(i,j)^{\rm th}$ entry is
$\Phi^{\D}(x_{ij}(\al))$.

\example\label{ex2} Let $\Q$ be
$$\vcenter
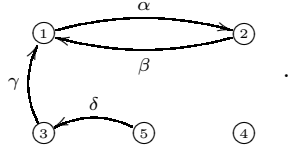
{
\xymatrix@C=1cm@R=1cm{ %
\vtx{1}\ar@/^/@{->}[rr]^{\al} \ar@/_/@{<-}[rr]_{\be}&&\vtx{2}\\
\vtx{3}\ar@/^/@{<-}[r]^{\de} \ar@/^/@{->}[u]^{\ga}&\vtx{5}& \vtx{4}\\
}} \quad.
$$
Define a mixed quiver setting $\QS=(\Q,\n,\g,\h,\i)$ by $\i(1)=2$, $\i(3)=4$, $\i(5)=5$;
$\g_1=\g_2=GL$, $\g_3=\g_4=SL$, $\g_5=O$; $\h_\al=\h_\ga=\h_\de=M$, $\h_\be=S^{+}$. Then
$\Q^{\D}$ is
$$
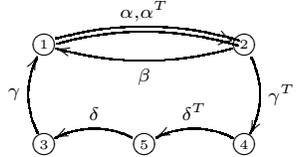
\vcenter{
\xymatrix@C=1cm@R=1cm{ %
\vtx{1}\ar@2@/^/@{->}[rr]^{\al,\al^T} \ar@/_/@{<-}[rr]_{\be}&&\vtx{2}\\
\vtx{3}\ar@/^/@{<-}[r]^{\de} \ar@/^/@{->}[u]^{\ga} &\vtx{5}&
\vtx{4}\ar@/_/@{->}[l]_{\de^T} \ar@/_/@{<-}[u]_{\ga^T} \\
}} .
$$
\bigskip %

The following theorem was proven by Zubkov in~\cite{ZubkovI}.

\begin{theo}\label{theo_ZubI} Let $(\Q,\n,\g,\h,\i)$ be a mixed quiver setting
satisfying~\Ref{eq_condition}. If $\g_v\in\{GL, O, \Sp\}$ for all $v\in\Q_0$, then
$\FF$-algebra $\FF[H(\Q,\n,\h)]^{G(\n,\g,\i)}$ is
generated by %
$$\Phi^{\D}(\si_t(X_{\be_1}\cdots X_{\be_s})),$$ %
where $\be_1\cdots \be_s$ is a closed path in $\Q^{\D}$ and $1\leq t\leq \n_{\be'_1}$.
\end{theo}

\remark\label{rem1} For any system of generators from Theorem~\ref{theo_ZubI} the exists
a {\it finite} subset that generates the corresponding algebra of invariants. %
\bigskip

Note that Theorem~\ref{theo_ZubI} implies all parts of Theorem~\ref{theo_matrix} but
part~d).

\subsection{Relations for invariants of mixed representations of quivers}\label{section_mixed_quivers_relations}

In this section we assume that $\QS=(\Q,\n,\g,\h,\i)$ is a mixed quiver setting with %
$$\g_v=GL\text{ for all }v\in\Q_0 \text{ and } \h_{\al}=M \text{ for all }\al\in\Q_1.$$ %
Then $\QS$-mixed
representations are called {\it mixed} representations (see Example~\ref{ex3part4}). For
short, we denote $H=H(\Q,\n,\h)$ and $G=GL(\n,\g,\i)$.


Denote by $\M$ the set of closed paths in the double quiver $\Q^{\D}$ (see Definition~\ref{def_double}).
Assume $(\al,\be,\ga)$ is a triple of paths in $\Q^{\D}$; then
it is called {\it incident} to $v\in\Q_0$ if $v<\i(v)$, $\al$ is incident to $v$, $\be$
is a path from $\i(v)$ to $v$, and $\ga$ is a path from $v$ to $\i(v)$. In what follows we will also need the notion of primitive path that was introduced in Section~\ref{subsection_def_quivers}.   

The following lemma is an analogue of Lemma~\ref{lemma_primitive}.

\begin{lemma}\label{lemma_primitive3}
\begin{enumerate}
\item[a)] The algebra $\FF[H]^{G}$ is generated by $\si_t(X_{\al})$, where $1\leq t\leq
\n_{\be'_1}$ and $\al\in\M$ is primitive.

\item[b)] If $a_i\in\FF$, $\al_i\in\M$ is incident to $v\in\Q_0$, and $t>0$, then the
element $\si_t(\sum_i a_i X_{\al_i})$ is a polynomial in $\si_l(X_{\al})$, where $l$ ranges over $\{1,\ldots,\n_v\}$ and $\al$ ranges over primitive elements from $\M$.

\item[c)] Let  $a_i,b_j,c_k\in\FF$, $(\al_i,\be_j,\ga_k)$ be a triple incident to $v\in\Q_0$,
and $t,r\geq0$. Then the element $\si_{t,r}(\sum_i a_i X_{\al_i},\,\sum_j b_j X_{\be_j},\,\sum_k c_k X_{\ga_k})$ is a polynomial in $\si_l(X_{\al})$, where $l$ ranges over $\{1,\ldots,\n_v\}$ and $\al$
ranges over primitive elements from $\M$.
\end{enumerate}
\end{lemma}
\bigskip %

Given $v\in\Q_0$, we consider the following equalities:
\begin{enumerate}
\item[a)] $\si_t(X_{\al}X_{\be})-\si_t(X_{\be}X_{\al})=0$ for $1\leq t\leq \n_v$, where $\al,\be$ are paths in $\Q^{\D}$ such that $\al\be\in\M$ is incident to $v$;

\item[b)] $\si_t(\sum_i a_i X_{\al_i})=0$ for $t>\n_v$, where $\al_i\in\M$ is incident to $v$;

\item[c)] $\si_{t,r}(\sum_i a_i X_{\al_i},\,\sum_j b_j X_{\be_j},\,\sum_k c_k X_{\ga_k})=0$ for $t+2r>n$, where $v<\i(v)$ and $(\al_i,\be_j,\ga_k)$ is a triple incident to $v$ for all $i,j,k$.

\item[d)] free relations, i.e., relations between generators of $\FF[H]^G$ that are valid for any dimension vector $\n$.
\end{enumerate}
Here we assume that Lemma~\ref{lemma_primitive3} is applied to these equalities, i.e.,
the left hand sides of the above equalities are polynomials in $\si_l(A)$, where $l$
ranges over $\{1,\ldots,\n_v\}$ and $\al$ ranges over primitive elements from $\M$.

\begin{theo}\label{theo_mixed quivers_relations} Consider the system of generators for $\FF[H]^G$
from part~a) of Lemma~\ref{lemma_primitive3}. Then the ideal of relations between them is
generated by relations a), b), c), and d) considered for all vertices of $\Q$.
\end{theo}

\remark{} Theorem~\ref{theo_mixed quivers_relations} was proved by Zubkov
in~\cite{ZubkovII}.

\section{Tableaux with substitutions and $\F$}\label{section_bpf}

\definition{of a block partial linearization of the pfaffian} Fix
$\un{n}=(n_1,\ldots,n_m)\in\NN^{m}$, where $n=n_1+\cdots+n_m$ is even. For any $1\leq
p,q\leq m$ and an $n_p\times n_q$ matrix $X$ denote by $X^{p,q}$ the $n\times n$ matrix,
partitioned into $m\times m$ number of blocks, where the block in the $(i,j)^{\rm th}$
position is an $n_i\times n_j$ matrix; the block in the $(p,q)^{\rm th}$ position is
equal to $X$, and the rest of blocks are zero matrices.

Let $1\leq p_1,\ldots,p_s,q_1,\ldots,q_s\leq m$, let $X_j$ be an $n_{p_j}\times n_{q_j}$
matrix for any $1\leq j\leq s$, and let $r_1,\ldots,r_s$ be positive integers, satisfying
$r_1+\cdots+r_s=n/2$. The element %
\begin{eq}\label{eq_block_linearization}
{\P}_{r_1,\ldots,r_s}(X_1^{p_1,q_1},\ldots,X_s^{p_s,q_s})
\end{eq}%
is a partial linearization of the pfaffian of block matrices
$X_1^{p_1,q_1},\ldots,X_s^{p_s,q_s}\!\!$, and it is  called a {\it block partial
linearization of the pfaffian} ({\it b.p.l.p.}).
\bigskip

We consider only b.p.l.p.-s, satisfying %
\begin{eq}\label{eq_block_condition}
\sum_{1\leq j\leq s,\;p_j=i}r_j + %
\sum_{1\leq j\leq s,\;q_j=i}r_j=n_i\;{\rm for\;all}\; 1\leq i\leq m,
\end{eq}%
since only these elements appear in the context of the invariant theory.

\example\label{ex pf_0} Let $\un{n}=(5,3)$ and let $X_1,X_2,X_3$
be $5\times 3$, $5\times 5$, and $3\times 3$ matrices respectively. Then %
$$X_1^{1,2}=\left(
\begin{array}{cc}
0&X_1\\
0&0\\
\end{array}
\right),\quad %
X_2^{1,1}=\left(
\begin{array}{cc}
X_2&0\\
0&0\\
\end{array}
\right),\quad %
X_3^{2,2}=\left(
\begin{array}{cc}
0&0\\
0&X_3\\
\end{array}
\right)
$$
are $8\times 8$ matrices. The b.p.l.p.~$f={\P}_{1,2,1}(X_1^{1,2},X_2^{1,1},X_3^{2,2})$
satisfies condition~\Ref{eq_block_condition}. It is convenient to display the data that
determine $f$ as a two-column tableau filled with
arrows: %
$$
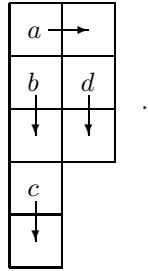
\begin{picture}(0,100)
\put(-10,90){%
\put(0,0){\rectangle{10}{10}\put(-3,-3){$a$}\put(5,0){\vector(1,0){15}}}%
\put(20,0){\rectangle{10}{10}\put(-3,-3){}}%
\put(0,-20){\rectangle{10}{10}\put(-3,-3){$b$}\put(0,-5){\vector(0,-1){15}}}%
\put(20,-20){\rectangle{10}{10}\put(-3,-3){$d$}\put(0,-5){\vector(0,-1){15}}}%
\put(0,-40){\rectangle{10}{10}\put(-3,-3){}}%
\put(20,-40){\rectangle{10}{10}\put(-3,-3){}}%
\put(0,-60){\rectangle{10}{10}\put(-3,-3){$c$}\put(0,-5){\vector(0,-1){15}}}%
\put(0,-80){\rectangle{10}{10}\put(-3,-3){}}%
\put(40,-30){.}
}%
\end{picture}$$%
\noindent The arrow $a$ determines the block matrix $X_1^{1,2}$ as follows: $a$ goes from
the $1^{\rm st}$ column to the $2^{\rm nd}$ column and we assign the matrix $X_1$ to $a$.
To determine the block matrix $X_2^{1,1}$ we take {\it two} arrows $b,c$ that go from the
$1^{\rm st}$ column to the $1^{\rm st}$ column, since the degree of $f$ in entries of
$X_2$ is $2$. Assign the matrix $X_2$ to $b$ and to $c$. Finally, the arrow $d$
determines the block matrix $X_2^{2,2}$: $d$ goes from the $2^{\rm nd}$ column to the
$2^{\rm nd}$ column and $X_2$ is assigned to it. Condition~\Ref{eq_block_condition}
implies that for every column the total amount of arrows that start or terminate in the
column is equal to the length of the column. Note that this conditions do not uniquely
determine arrows of a tableau.
\bigskip

Below we formulate the definition of a tableau with substitution that gives alternative
way to work with b.p.l.p.-s. Acting in the same way as in Example~\ref{ex pf_0}, for
every b.p.l.p.~$f$ given by~\Ref{eq_block_linearization} we construct a tableau with
substitution $(\T,(X_1,\ldots,X_s))$ and define the function ${\F}_{\T}$ such that $f=\pm
{\F}_{\T}(X_1,\ldots,X_s)$.

\definition{of shapes} The {\it shape} of dimension
$\un{n}=(n_1,\ldots,n_m)\in\NN^{m}$ is the collection of $m$ columns of cells. The
columns are numbered by $1,2,\ldots,m$, and the $i^{\rm th}$ column contains exactly
$n_i$ cells, where $1\leq i\leq m$. Numbers $1,\ldots,n_i$ are assigned to the cells of
the $i^{\rm th}$ column, starting from the top. As an example, the shape of
dimension $\un{n}=(3,2,3,1,1)$ is%
$$\begin{picture}(0,70)%
\put(-40,45){%
\put(0,0){\rectangle{10}{10}\put(-2,-4){1}}%
\put(20,0){\rectangle{10}{10}\put(-2,-4){1}}%
\put(40,0){\rectangle{10}{10}\put(-2,-4){1}}%
\put(60,0){\rectangle{10}{10}\put(-2,-4){1}}%
\put(80,0){\rectangle{10}{10}\put(-2,-4){1}}%
\put(0,-20){\rectangle{10}{10}\put(-3,-4){2}}%
\put(20,-20){\rectangle{10}{10}\put(-3,-4){2}}%
\put(40,-20){\rectangle{10}{10}\put(-3,-4){2}}%
\put(0,-40){\rectangle{10}{10}\put(-3,-4){3}}%
\put(40,-40){\rectangle{10}{10}\put(-3,-4){3}}%
\put(-2,17){1}%
\put(17,17){2}%
\put(37,17){3}%
\put(57,17){4}%
\put(77,17){5}%
}%
\end{picture}$$ %

\definition{of a tableau with substitution} Let
$\un{n}=(n_1,\ldots,n_m)\in\NN^{m}$ and let $n=n_1+\cdots+n_m$ be even. A pair
$(\T,(X_1,\ldots,X_s))$ is called a {\it tableau with substitution} of dimension $\un{n}$
if
\begin{enumerate}
\item[$\bullet$] $\T$ is the shape of dimension $\un{n}$ together with a set of arrows. An
{\it arrow} goes from one cell of the shape into another one, and each cell of the shape
is either the head or the tail of one and only one arrow.  We refer to $\T$ as {\it
tableau} of dimension $\un{n}$, and we write $a\in \T$ for an arrow $a$ from $\T$. Given an
arrow $a\in \T$, denote by $a'$ and $a''$ the columns containing the head and the tail of
$a$, respectively. Similarly, denote by $'a$ the number assigned to the cell containing
the head of $a$, and denote by $''a$ the number assigned to the cell containing the tail
of $a$. Schematically this is depicted as
$$
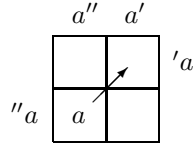
\begin{picture}(0,50)
\put(-10,30){%
\put(0,0){\rectangle{10}{10}\put(-3,-3){}}%
\put(20,0){\rectangle{10}{10}\put(-3,-3){}}%
\put(0,-20){\rectangle{10}{10}\put(-3,-3){$a$}\put(5,5){\vector(1,1){13}}}%
\put(20,-20){\rectangle{10}{10}\put(-3,-3){}}%
\put(-3,15){$a''$}%
\put(17,15){$a'$}%
\put(35,-3){$'a$}%
\put(-26,-23){$''a$}%
}%
\end{picture}
$$%

\item[$\bullet$] $\varphi$ is a fixed mapping from the set of arrows of $\T$ onto $[1,s]$
that satisfies the following property:

\begin{enumerate}
\item[] if $a,b\in \T$ and $\ovphi{a}=\ovphi{b}$, then $a'=b'$, $a''=b''$;
\end{enumerate}

\item[$\bullet$] $(X_1,\ldots,X_s)$ is a sequence of matrices such that the matrix
$X_{\ovphi{a}}$ assigned to the arrow $a\in \T$ is $n_{a''}\times n_{a'}$ matrix and its
$(p,q)^{\rm th}$ entry is denoted by $(X_j)_{pq}$.
\end{enumerate}

\example\label{ex_bp_new} Let $\T$ be the tableau
$$
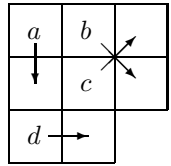
\begin{picture}(0,65)
\put(-20,50){%
\put(0,0){\rectangle{10}{10}\put(-3,-3){$a$}\put(0,-5){\vector(0,-1){15}}}%
\put(20,0){\rectangle{10}{10}\put(-3,-3){$b$}\put(5,-5){\vector(1,-1){13}}}%
\put(40,0){\rectangle{10}{10}}%
\put(0,-20){\rectangle{10}{10}}%
\put(20,-20){\rectangle{10}{10}\put(-3,-3){$c$}\put(5,5){\vector(1,1){13}}}%
\put(40,-20){\rectangle{10}{10}}%
\put(0,-40){\rectangle{10}{10}\put(-3,-3){$d$}\put(5,0){\vector(1,0){15}}}%
\put(20,-40){\rectangle{10}{10}}%
}%
\end{picture}
$$
of dimension $(3,3,2)$. Define $\varphi$ by $\ovphi{a}=1$, $\ovphi{b}=\ovphi{c}=2$, and
$\ovphi{d}=3$, and let $X_1$, $X_3$ be $3\times 3$ matrices and $X_2$ be a $3\times 2$
matrix. Then $(\T,(X_1,X_2,X_3))$ is a tableau with substitution.

\definition{of ${\F}_{\T}(X_1,\ldots,X_s))$} Let $(\T,(X_1,\ldots,X_s))$ be
a tableau with substitution of dimension
$\un{n}$. Define the polynomial %
\begin{eq}\label{eq_bp_def_F0}
{\F}_{\T}^0(X_1,\ldots,X_s)=\sum_{\pi_1\in \Symmgr_{n_1},\ldots,\pi_m\in \Symmgr_{n_m}}
\sign(\pi_1)\cdots\sign(\pi_m)\prod_{a\in \T}
(X_{\ovphi{a}})_{\pi_{a''}(''a),\pi_{a'}('a)},
\end{eq}%
and the coefficient
$$c_{\T}=\prod_{j=1}^s\#\{a\in \T\,|\,\ovphi{a}=j\}!$$ %
In the case $\FF=\QQ$ define
$${\F}_{\T}(X_1,\ldots,X_s)=\frac{1}{c_{\T}}{\F}^0_{\T}(X_1,\ldots,X_s).$$ %
Since ${\F}_{\T}(X_1,\ldots,X_s)$ is a polynomial in entries of $X_1,\ldots,X_s$ with
integer coefficients, the definition of ${\F}_{\T}(X_1,\ldots,X_s)$ extends over an
arbitrary field.
\bigskip

\noindent{}The next lemma shows that $\F$ is a b.p.l.p.

\begin{lemma}\label{lemma_bplp}
\begin{enumerate}
\item[a)] Let $f$ be a b.p.l.p.~\Ref{eq_block_linearization}
satisfying~\Ref{eq_block_condition}. Then there is a tableau with substitution
$(\T,(X_1,\ldots,X_s))$ such that ${\F}_{\T}(X_1,\ldots,X_s)=\pm f$.

\item[b)] For every tableau with substitution $(\T,(X_1,\ldots,X_s))$ of dimension
$\un{n}$ there is a b.p.l.p.~$f$ satisfying~\Ref{eq_block_condition} such that
${\F}_{\T}(X_1,\ldots,X_s)=\pm f$.
\end{enumerate}
\end{lemma}

\begin{proof}
{\bf a)} A tableau with substitution $(\T,(X_1,\ldots,X_s))$ of dimension $\un{n}$ is
constructed as follows. Its arrows are $a_{j r}$, where $a_{j r}$ goes from the $p_j^{\rm
th}$ to the $q_j^{\rm th}$ column for $1\leq j\leq s$ and $1\leq r\leq k_j$.
Condition~\Ref{eq_block_condition} guarantees that for any $1\leq i\leq m$ the total
number of arrows that begin or end in the $i^{\rm th}$ column is $n_i$. Complete the
construction by setting $\ovphi{a_{j r}}=j$. Note that $(\T,(X_1,\ldots,X_s))$ is not
uniquely determined by $f$.

Since $X_1^{p_1,q_1},\ldots,X_s^{p_s,q_s}$ are block matrices, the permutation $\pi$ from
formula~\Ref{eq_lin_P} can be written as the composition
$\pi=\pi_1\circ\cdots\circ\pi_m$, where the permutation $\pi_i$ acts as an identity on
the set $[1,c_{i-1}]\,\cup\,[c_{i+1}+1,n/2]$ for $c_i=k_1+\cdots+k_i$ ($1\leq i\leq m$).
The claim follows from formula~\Ref{eq_lin_P}.

{\bf b)} Consider $a_1,\ldots,a_s\in \T$ such that $\ovphi{a_1}=1,\ldots,\ovphi{a_s}=s$.
Then
$$\F_{\T}(X_1,\ldots,X_s)=\pm\P_{r_1,\ldots,r_s}(X_1^{a_1'',a_1'},\ldots,X_s^{a_s'',a_s'}),$$
where $r_j=\#\{a\in \T\,|\,\ovphi{a}=j\}$ for any $1\leq j\leq s$.
\end{proof}

\example\label{ex_pf_1} {\bf 1.} Let $n$ be even and $(\T,(X_1,\ldots,X_s))$ be a tableau
with substitution of dimension $(n)$, where $\T$ is
$$
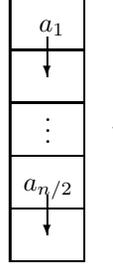
\begin{picture}(0,100)
\put(0,90){%
\put(0,0){\rectangle{14}{10}\put(-3,-3){$a_1$}\put(0,-5){\vector(0,-1){15}}}%
\put(0,-20){\rectangle{14}{10}\put(-3,-3){}}%
\put(0,-40){\rectangle{14}{10}\put(-1,-5){$\vdots$}}%
\put(0,-60){\rectangle{14}{10}\put(-9,-3){$a_{n/2}$}\put(0,-5){\vector(0,-1){15}}}%
\put(0,-80){\rectangle{14}{10}\put(-3,-3){}}%
\put(24,-40){.}
}%
\end{picture}
$$
In other words, arrows of $\T$ are $a_1,\ldots,a_{n/2}$, where $'a_i=2i$, $''a_i=2i-1$,
$a_i'=a_i''=1$ for $1\leq i\leq n/2$, and $X_1,\ldots,X_s$ are $n\times n$ matrices.

If $s=1$, then~\Ref{eq_P} implies $\F_{\T}(X_1)=\P(X_1)$ for $\FF=\QQ$ and consequently for
an arbitrary $\FF$.

If $s>1$, then by~\Ref{eq_lin_P} %
$${\F}_{\T}(X_1,\ldots,X_s)={\P}_{r_1,\ldots,r_s}(X_1,\ldots,X_s),$$ %
where $r_j=\#\{a\in \T\,|\,\ovphi{a}=j\}$ for any $1\leq j\leq s$, is a partial
linearization of the pfaffian.
\smallskip

{\bf 2.} For $n\times n$ matrices $X_1,\ldots,X_s$ let $(\T,(X_1,\ldots,X_s))$ be the
tableau with substitution of dimension $(n,n)$, where $\T$ is
$$
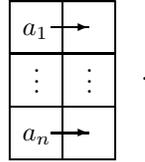
\begin{picture}(0,50)
\put(-10,40){%
\put(0,0){\rectangle{10}{10}\put(-5,-3){$a_1$}\put(6,0){\vector(1,0){14}}}%
\put(20,0){\rectangle{10}{10}\put(-3,-3){}}%
\put(0,-20){\rectangle{10}{10}\put(-1,-5){$\vdots$}}%
\put(20,-20){\rectangle{10}{10}\put(-1,-5){$\vdots$}}%
\put(0,-40){\rectangle{10}{10}\put(-5,-3){$a_n$}\put(6,0){\vector(1,0){14}}}%
\put(20,-40){\rectangle{10}{10}\put(-3,-3){}}%
\put(40,-20){.}
}%
\end{picture}
$$
If $s=1$, then the formula
\begin{eq}\label{eq_det}
\det(X)=\frac{1}{n!}\sum_{\pi_1,\pi_2\in \Symmgr_n}\sign(\pi_1)\sign(\pi_2) \prod_{i=1}^n
x_{\pi_1(i),\pi_2(i)},
\end{eq}
which is valid over $\QQ$, implies the equality $\F_{\T}(X_1)=\det(X_1)$ over every $\FF$.
For $s>1$ the expression $\F_{\T}(X_1,\ldots,X_s)$ is a partial linearization of the
determinant.
\smallskip

{\bf 3.} If $(\T,(X_1,X_2))$ is a tableau with substitution from part~2 of
Example~\ref{ex_pf_1}, where for $1\leq k\leq n$ we have
$\ovphi{a_1}=\cdots=\ovphi{a_k}=1$, $\ovphi{a_{k+1}}=\cdots=\ovphi{a_n}=2$, $X_1=X$, and
$X_2=E$ is the identity $n\times n$ matrix, then ${\F}_{\T}(X,E)=\si_k(X)$.
\smallskip %

{\bf 4.} Suppose that $t,r,s\in\NN$, and $X$, $Y$, $Z$, respectively, are matrices of
dimensions $(t+2r)\times(t+2s)$, $(t+2r)\times(t+2r)$, $(t+2s)\times(t+2s)$,
respectively. Let $\T$ be a tableau of dimension $(t+2r,t+2s)$ that in the case
$r=s$ is
depicted as %
$$\begin{picture}(0,160)
\put(-10,150){%
\put(0,0){\rectangle{10}{10}\put(-5,-3){$a_1$}\put(6,0){\vector(1,0){14}}}%
\put(20,0){\rectangle{10}{10}}%
\put(0,-20){\rectangle{10}{10}\put(-1,-5){$\vdots$}}%
\put(20,-20){\rectangle{10}{10}\put(-1,-5){$\vdots$}}%
\put(0,-40){\rectangle{10}{10}\put(-5,-3){$a_t$}\put(6,0){\vector(1,0){14}}}%
\put(20,-40){\rectangle{10}{10}}%
\put(0,-60){\rectangle{10}{10}\put(-4,-2){$b_1$}\put(0,-5){\vector(0,-1){15}}}%
\put(0,-80){\rectangle{10}{10}}%
\put(0,-100){\rectangle{10}{10}\put(-1,-5){$\vdots$}}%
\put(0,-120){\rectangle{10}{10}\put(-4,-2){$b_r$}\put(0,-5){\vector(0,-1){15}}}%
\put(0,-140){\rectangle{10}{10}}%
\put(20,-60){\rectangle{10}{10}\put(-4,-2){$c_1$}\put(0,-5){\vector(0,-1){15}}}%
\put(20,-80){\rectangle{10}{10}}%
\put(20,-100){\rectangle{10}{10}\put(-1,-5){$\vdots$}}%
\put(20,-120){\rectangle{10}{10}\put(-4,-2){$c_s$}\put(0,-5){\vector(0,-1){15}}}%
\put(20,-140){\rectangle{10}{10}}%
\put(40,-70){,}
}%
\end{picture}
$$
otherwise define $\T$ analogously. Let $(\T,(X,Y,Z))$ be a tableau with substitution, where
$\ovphi{a_i}=1$, $\ovphi{b_j}=2$, and $\ovphi{c_k}=3$ for $1\leq i\leq t$, $1\leq j\leq
r$, and $1\leq k\leq s$. Then $\F_{\T}(X,Y,Z)=q\,\DP_{r,s}(X,Y,Z)$, where $q=\pm1$
and $\DP_{r,s}$ was introduced in Section~3 of~\cite{LZ1}. Moreover, if $r=s$,
then $q=1$.

The above mentioned polynomial $\DP_{r,s}$ is called {\it determinant-pfaffian}. This function inherits some properties from the determinant as well as from the pfaffian:
\begin{enumerate}
\item[a)] $\DP_{0,0}(X,Y,Z)=\det(X)$;

\item[b)] $\DP_{r,s}(X,Y,Z)=\P(Y)\,\P(Z)$ if $t=0$;

\item[${\rm c')}$]   $\DP_{r,s}(gX,gYg^T,Z)=\det(g)\, \DP_{r,s}(X,Y,Z)$, where $g$ is a   $(t+2r)\times(t+2r)$ matrix;  

\item[${\rm c'')}$]  $\DP_{r,s}(Xg,Y,g^TZg)=\det(g)\, \DP_{r,s}(X,Y,Z)$, where $g$ is a  $(t+2s)\times(t+2s)$ matrix. 
\end{enumerate}
\bigskip

The following formula is the main result of~\cite{Lop_bplp} and it plays a crucial role
in proofs of results on $\F$ and invariants.

\begin{theo}\label{theo_short}
{\rm (Decomposition formula: short version).}\\
Let $(\T,(X_1,\ldots,X_s))$ be a tableau with substitution of dimension $\un{n}\in\NN^m$.
Let $1\leq q_1<q_2\leq m$, $n_{q_1}=n_{q_2}$, and the vector $\un{d}\in\NN^{m-2}$ be
obtained from $\un{n}$ by eliminating the $q_1^{\it th}$ and the $q_2^{\it th}$
coordinates. Then $\F_{\T}(X_1,\ldots,X_s)$ is a polynomial in $\F_{\mathcal D}(Y_1,\ldots,Y_l)$ and $\si_t(h)$, where
\begin{enumerate}
\item[$\bullet$] $({\mathcal D},(Y_1,\ldots,Y_l))$ ranges over tableaux with substitutions of dimension $\un{d}$ such that $Y_1,\ldots,Y_l$ are products of matrices
$X_1,\ldots,X_s,X_1^T,\ldots,X_s^T$;

\item[$\bullet$] $h$ ranges over products of matrices
$X_1,\ldots,X_s,X_1^T,\ldots,X_s^T$;

\item[$\bullet$] $t$ ranges over $[1,n_{q_1}]$. %
\end{enumerate} Moreover,
coefficients of this polynomial belong to the image of $\ZZ$ in $\FF$ under the natural
homomorphism.
\end{theo}
\bigskip %

\noindent The explicit formulation of the decomposition formula can be found
in~\cite{Lop_bplp}.

\section{Generators for invariants of $\QS$-mixed representations of quiver}\label{section_generators}

In this section we describe generators for invariants of arbitrary mixed quiver setting $\QS$.
We start with a particular case of semi-invariants for bipartite quivers.

\subsection{Generators for semi-invariants of representations of bipartite quivers}\label{subsection_semi}

Given a quiver setting $(\Q,\n)$, denote
$$SL(\n)=\prod_{v\in \Q_0} SL(\n_v).$$ %
Then $\FF[H(\Q,\n)]^{SL(\n)}$ is the algebra of {\it semi-invariants} of $\Q$.

A quiver $\Q$ is called {\it bipartite}, if every vertex is either a source (i.e.~there
is no arrow ending at this vertex) or a sink (i.e.~there is no arrow starting at this
vertex).

\definition{of a $(\Q,\n)$-tableau with substitution} A tableau with
substitution $(\T,(Y_1,\ldots,Y_s))$ of dimension $\un{n}\in\NN^m$ is called a {\it
$(\Q,\n)$-tableau with substitution}, if for some {\it weight}
$\un{w}=(w_1,\ldots,w_l)\in\NN^l$ we have
\begin{enumerate}
\item[$\bullet$] $\T$ is a union of $m$ rectangular blocks $B_1,\ldots,B_l$, where $B_i$
consists of $w_i$ columns of length $\n_i$ ($1\leq i\leq l$); this condition is
equivalent to
$\un{n}=(\underbrace{\n_1,\ldots,\n_1}_{w_1},\ldots,\underbrace{\n_l\ldots,\n_l}_{w_l})$;

\item[$\bullet$] if $a\in \T$, then there exists an $\al\in \Q_1$ such that
\begin{enumerate}
\item[a)] $X_\al=Y_{\ovphi{a}}$,

\item[c)] $B_{\al'}$ contains column $a''$ of $\T$,

\item[b)] $B_{\al''}$ contains column $a'$ of $\T$.
\end{enumerate}
\end{enumerate}

\begin{theo}\label{theo_zigzag}
Let $(\Q,\n)$ be a quiver setting such that $\Q$ is bipartite. Then the algebra of
semi-invariants $\FF[H(\Q,\n)]^{SL(\n)}$ is spanned over $\FF$ by the elements
$\F_D(Z_1,\ldots,Z_h)$, where $(D,(Z_1,\ldots,Z_h))$ is a $(\Q,\n)$-tableau with
substitution.
\end{theo}
\bigskip %
Let us remark that under the given restrictions on the quiver the generating set from
Theorem~\ref{theo_zigzag} is smaller than that from Theorem~\ref{theo_main} (see below).

\remark\label{rem_history4} Generators for semi-invariants of an {\it arbitrary} quiver
were established by Domokos and Zubkov in~\cite{DZ01} using the methods
from~\cite{Donkin92a},~\cite{Donkin94},~\cite{Zubkov96},~\cite{Zubkov_Fund_Math_01}, and,
independently, by Derksen and Weyman in~\cite{DW_LR_02},~\cite{DW00} utilizing the
methods of the representation theory of quivers. Simultaneously, similar result in the
case of characteristic zero was obtained by Schofield and van den Bergh
in~\cite{Schofield_van_den_Bergh_01}.

\example\label{ex2_zigzag} Consider the following quiver $\Q$:
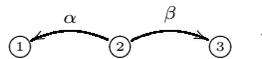
$$\vcenter{
\xymatrix@C=1cm@R=1cm{ %
\vtx{1}\ar@/^/@{<-}[r]^{\al} & \vtx{2} \ar@/^/@{->}[r]^{\be}&\vtx{3}\\
}} \quad.
$$
Let $(\T,(Y_1,\ldots,Y_s))$ be a $(\Q,\n)$-tableau with substitution of a weight
$\un{w}=(w_1,w_2,w_3)$. Then $\T$ consists of three blocks $B_1,B_2,B_3$, where $B_i$ is
the $\n_i\times w_i$ rectangle ($1\leq i\leq3$). Schematically $\T$ is depicted as
$$
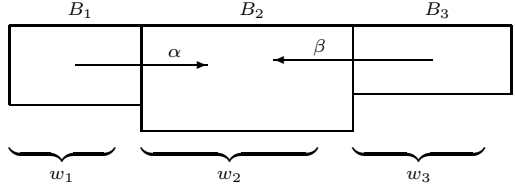
\begin{picture}(0,80)
\put(0,35){%
\put(-65,5){\rectangle{25}{15}}%
\put(0,0){\rectangle{40}{20}}%
\put(70,7){\rectangle{30}{13}}%
\put(-68,24){$\scriptstyle B_1$}%
\put(-3,24){$\scriptstyle B_2$}%
\put(67,24){$\scriptstyle B_3$}%
\put(-65,5){\vector(1,0){50}}%
\put(70,7){\vector(-1,0){60}}%
\put(-30,8){$\scriptstyle \al$}%
\put(25,10){$\scriptstyle \be$}%
\put(-40,-25){$\underbrace{\qquad\qquad\qquad\quad\!\!}_{w_2}$}%
\put(-90,-25){$\underbrace{\qquad\qquad}_{w_1}$}%
\put(40,-25){$\underbrace{\qquad\qquad\quad}_{w_3}$}%
\put(115,0){.}%
}%
\end{picture}
$$
Here arrows of $\T$ are denoted by the same letters as the corresponding arrows of $\Q$.

Now let us to consider a concrete example. We assume that $\n=(1,2,1)$ and define a
tableau with substitution $(\T,(Y_1,\ldots,Y_4))$ of dimension $\un{n}=(1,2,2,1,1,1)$ by
$\T$:
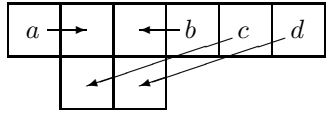
$$\begin{picture}(0,45)
\put(-50,30){%
\put(0,0){\rectangle{10}{10}\put(-3,-3){$a$}\put(5,0){\vector(1,0){15}}}%
\put(20,0){\rectangle{10}{10}}%
\put(20,-20){\rectangle{10}{10}}%
\put(40,0){\rectangle{10}{10}}%
\put(40,-20){\rectangle{10}{10}}%
\put(60,0){\rectangle{10}{10}\put(-3,-3){$b$}\put(-5,0){\vector(-1,0){15}}}%
\put(80,0){\rectangle{10}{10}\put(-3,-3){$c$}\put(-5,-3){\vector(-3,-1){55}}}%
\put(100,0){\rectangle{10}{10}\put(-3,-3){$d$}\put(-5,-3){\vector(-3,-1){55}}}%
}%
\end{picture}
$$
and by equalities $Y_{\varphi(a)}=X_{\al}$, %
$Y_{\varphi(b)}=Y_{\varphi(c)}=Y_{\varphi(d)}=X_{\be}$. Then $(\T,(Y_1,\ldots,Y_4))$ is
$(\Q,\n)$-tableau with substitution of the weight $\un{w}=(1,2,3)$.

\subsection{General case}\label{subsection_general_case}
Let $\QS=(\Q,\n,\g,\h,\i)$ be a mixed quiver setting and $\Q_0=\{1,\ldots,l\}$. As in
Section~\ref{subsection_mixed_quivers_generators}, without loss of generality we can
assume that condition~\Ref{eq_condition} is valid. In this section we will use the notion of mixed double quiver setting $\QS^{\D}$ (see Definition~\ref{def_double}). 

\definition{of a path $\QS$-tableau with substitution} A tableau with
substitution $(\T,(Y_1,\ldots,Y_s))$ of dimension $\un{n}\in\NN^m$ is called a {\it path
$\QS$-tableau with substitution}, if for some {\it weight}
$\un{w}=(w_1,\ldots,w_l)\in\NN^l$  we have
\begin{enumerate}
\item[$\bullet$] $\T$ is a union of $m$ rectangular blocks $B_1,\ldots,B_l$, where $B_i$
consists of $w_i$ columns of length $\n_i$ ($1\leq i\leq l$); this condition is
equivalent to
$\un{n}=(\underbrace{\n_1,\ldots,\n_1}_{w_1},\ldots,\underbrace{\n_l\ldots,\n_l}_{w_l})$;

\item[$\bullet$] if $a\in \T$, then there exists a path $\al=\al_1\cdots \al_r$ in $\Q$
(where $\al_1,\ldots,\al_r\in \Q_1$) such that
\begin{enumerate}
\item[a)] $X_{\al_1}\cdots X_{\al_r}=Y_{\ovphi{a}}$,

\item[b)] $B_{\al'}$ contains column $a''$ of $\T$,

\item[c)] $B_{\i(\al'')}$ contains column $a'$ of $\T$.
\end{enumerate}
\end{enumerate}

\example\label{ex2_1} Consider the mixed quiver setting $\QS$ defined in
Example~\ref{ex2}. Let $(\T,(Y_1,\ldots,Y_s))$ be a path $\QS^{\D}$-tableau with
substitution of a weight $\un{w}=(w_1,\ldots,w_5)$. Then $\T$ consists of five rectangular
blocks $B_1,\ldots,B_5$, where $\n_i\times w_i$ block $B_i$ corresponds to vertex $i$
($1\leq i\leq 5$). Schematically $\T$ is depicted as
$$
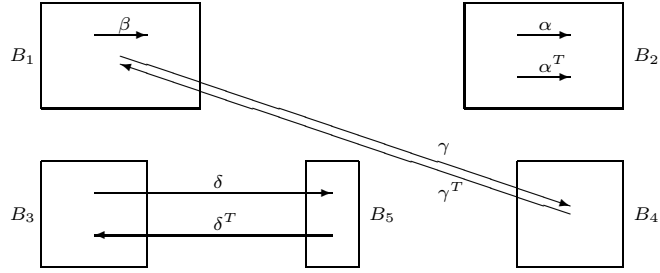
\begin{picture}(0,110)
\put(0,85){%
\put(-80,0){\put(0,0){\rectangle{30}{20}}%
\put(-10,8){\vector(1,0){20}}%
\put(-1,10){$\scriptstyle \be$}}%
\put(80,0){\put(0,0){\rectangle{30}{20}}%
\put(-10,8){\vector(1,0){20}}%
\put(-2,10){$\scriptstyle \al$}%
\put(-10,-8){\vector(1,0){20}}%
\put(-2,-6){$\scriptstyle \al^T$}}%
\put(-90,-60){\rectangle{20}{20}}%
\put(90,-60){\rectangle{20}{20}}%
\put(0,-60){\rectangle{10}{20}}%
\put(-80,0){\vector(3,-1){170}}%
\put(90,-60){\vector(-3,1){170}}%
\put(-90,-52){\vector(1,0){90}}\put(-45,-50){$\scriptstyle \de$}%
\put(0,-68){\vector(-1,0){90}}\put(-45,-66){$\scriptstyle \de^T$}%
\put(40,-36){$\scriptstyle \ga$}%
\put(40,-54){$\scriptstyle \ga^T$}%
\put(-122,-2){$\scriptstyle B_1$}%
\put(114,-2){$\scriptstyle B_2$}%
\put(-122,-62){$\scriptstyle B_3$}%
\put(114,-62){$\scriptstyle B_4$}%
\put(14,-62){$\scriptstyle B_5$}%
\put(130,-30){.}%
}%
\end{picture}
$$
Here we depicted arrows of $\T$ that correspond to paths in $\Q^{\D}$ of length one, i.e.,
arrows of $\Q^{\D}$; arrows of $\T$ are denoted by the same letters as the corresponding
arrows of $\Q^{\D}$.
\medskip %

\begin{theo}\label{theo_main} %
Let $(\Q,\n,\g,\h,\i)$ be a mixed quiver setting satisfying~\Ref{eq_condition}. Then the
algebra of invariants $\FF[H(\Q,\n,\h)]^{G(\n,\g,\i)}$ is generated as $\FF$-algebra by
the elements $\Phi^{\D}(\si_t(X_{\be_1}\cdots X_{\be_r}))$,
$\Phi^{\D}(\F_{\T}(Y_1,\ldots,Y_s))$, where
\begin{enumerate}
\item[1.] $\be_1\cdots \be_r$ ranges over all closed paths in $\Q^{\D}$ and
$1\leq t\leq \n_{\be'_1}$;    

\item[2.] $(\T,(Y_1,\ldots,Y_s))$ ranges over all path $\QS^{\D}$-tableaux with
substitutions of a weight $\un{w}$ such that 

\begin{enumerate} %
\item[a)] if $\g_v\in\{GL,O,\Sp\}$ for some $v\in \Q_0$, then
$w_{\i(v)}=w_{v}=0$;  

\item[b)] if $\g_v=SL$ for some $v\in \Q_0$, then $w_{\i(v)}=0$ or
$w_{v}=0$;            

\item[c)] if $\g_v=SO$ for some $v\in \Q_0$, then $w_{v}\leq 1$ and
$\i(v)=v$.            
\end{enumerate}
\end{enumerate}
\end{theo}
\bigskip %


\remark{} Semi-invariants of mixed representations of quivers were found by Lopatin and
Zubkov in~\cite{LZ1}. This result was generalized for an arbitrary quiver setting by
Lopatin in~\cite{Lop_so_inv}, where Theorem~\ref{theo_main} was proven. In particular,
Theorem~\ref{theo_main} implies the description of semi-invariants of {\it supermixed}
representations.

\bigskip
\noindent{\bf Acknowledgements.} This paper was supported by RFFI 07-01-00392.


\end{document}